\documentclass{amsart}

\usepackage{amsfonts}
\usepackage{amssymb}
\usepackage[forcekdg]{kuvio}

\newcommand{\N}{{\mathbb N}}
\newcommand{\Z}{{\mathbb Z}}

\newcommand{\A}{\mathcal{A}}
\newcommand{\ga}{$\Gamma(G,{\mathcal A})$}

\newcommand{\cN}{{\mathcal{N}}}
\newcommand{\R}{\mathcal{R}}
\newcommand{\F}{\mathcal{F}}

\theoremstyle{plain}
\newtheorem{lem}{Lemma}
\newtheorem{cor}{Corollary}
\newtheorem{thm}{Theorem}
\newtheorem{theorem}{Theorem}

\theoremstyle{remark}
\newtheorem{rem}{Remark}
\theoremstyle{definition}

\newcommand{\eqi}{\stackrel{ i}{=}}

\newcommand{\eqin}{\stackrel{ 0}{=}}

\begin{document}

\title{Periodic quotients of hyperbolic and large groups}

\author{A. Minasyan}
\address[Ashot Minasyan]{Universit\'{e} de Gen\`{e}ve,
Section de Math\'{e}matiques, 2-4 rue du Li\`{e}vre, Case postale
64, 1211 Gen\`{e}ve 4, Switzerland. On leave from School of
Mathematics, University of Southampton, Highfield, Southampton,
SO17 1BJ, UK.} \email{aminasyan@gmail.com}

%\thanks{The work of the first author was supported by the Swiss
%National Science Foundation Grants $\sharp$~PP002-68627, PP002-116899.}

\author{A.Yu. Olshanskii}
\address[Alexander Yu. Olshanskii]{Department of Mathematics,
1326 Stevenson Center, Vanderbilt University Nashville, TN 37240,
USA. And Department of Mechanics and Mathematics, Moscow State
University, Moscow, 119899 Russia.}
\email{alexander.olshanskiy@vanderbilt.edu}

\author{D. Sonkin}
\address[Dmitriy Sonkin]{Department of Mathematics, P. O. Box 400137
University of Virginia, Charlottesville, VA 22904, USA.}
\email{ds5nd@virginia.edu}

\thanks{The work of the first author was supported by the Swiss
National Science Foundation Grants $\sharp$~PP002-68627,
PP002-116899. The second author was supported by NSF Focus
Research Grant DMS-0455881, by NSF grant DMS-0700811 and by
Russian Fund for Fundamental Research grant $\sharp$~05-01-00895.}

\begin{abstract}  Let $G$ be either a non-elementary (word)
hyperbolic group or a large group (both in the sense of Gromov). In this paper we describe
several approaches for constructing continuous families of periodic
quotients of $G$ with various properties.

The first three methods work for any non-elementary hyperbolic
group, producing three different continua of periodic quotients of
$G$. They are based on the results and techniques, that were
developed by Ivanov and Olshanskii in order to show that there
exists an integer $n$ such that $G/G^n$ is an infinite group
of exponent $n$.

The fourth approach starts with a large group $G$ and produces a
continuum of pairwise non-isomorphic periodic residually finite
quotients.

Speaking of a particular application, we use each of these
methods to give a positive answer to a question of Wiegold from Kourovka Notebook. %\cite{Kourovka}.
\end{abstract}
%%2000 MSC:

\maketitle

%%%%%%%%%%%%%%%%%%%%%%%%%%%%%%%%%%%%%%%%%%%%%%%%%%%%%%%%%%%%%%%%%%%%%%%%%%%%%%%%%%%%%%%%%%%%%%%%%%%%%%%%%%%%%%%

\section{Introduction}
Let $G$ be an infinite finitely generated group. For $n \in \N$,
by $G^n$ we denote the (normal) subgroup generated by the $n$-th
powers of elements of $G$. Recall that $G$ is called {\it
periodic} if each element of $G$ has a finite order. The group $G$
is said to be of {\it bounded exponent} if there is $n \in N$ such
that $x^n=1$ for all $x \in G$, or, in other words, $G^n=\{1\}$;
the least integer $n$ for which this equality holds is called the
{\it exponent} of $G$. If for a periodic group $G$ such $n$ does
not exist, it is said to be {\it periodic of unbounded exponent}.
A group is called {\it elementary} if it has a cyclic subgroup of
finite index.

The question whether the quotient group $G/G^n$ (which is of
bounded exponent $n$) can be infinite has a long history in Group
Theory. In the case when $G=F_m$ is a non-abelian free group, it
was asked by W. Burnside in 1902 \cite{Burnside} and answered in
the affirmative much later by P. Novikov and S. Adian
\cite{Nov-Adian}
%\footnote{sm. takzhe knigu Adiana}
(see also \cite{Ad}) for any sufficiently large odd exponent $n$.
Yet later, A. Olshanskii \cite{Olsh-mono} suggested a different,
geometric approach to this problem, and, afterwards, in
\cite{Olsh1}, applied similar methods to show that for every
non-elementary torsion-free hyperbolic (in the sense of M. Gromov
\cite{Gromov}) group $G$ there exists an integer $n \in \N$ such
that $G/G^n$ is infinite. In \cite{Iv} S. Ivanov proved that
$F_m/F_m^n$ is infinite for any sufficiently large even exponent $n$
divisible by $2^9$ (infiniteness of groups $F_m/F_m^n$ for large
even values of $n$ was also proved by I. Lys\"enok \cite{Lys}).
Finally, in \cite{Iv-Olsh} Ivanov and Olshanskii combined the
techniques from \cite{Olsh1} and \cite{Iv} and obtained the
following theorem:

\begin{theorem} \label{thm:Iv-Olsh} For every non-elementary hyperbolic group G there exists a positive
even integer $n = n(G)$ such that the following are true:
\begin{itemize}
\item[(a)] The quotient group $G/G^n$ is infinite.
%\item[(b)] The word and conjugacy problems are solvable in $G/G^n$.
\item[(b)] Suppose $n = n_1n_2$, where $n_1$ is odd and $n_2$ is a
power of $2$. Then every finite subgroup of $G/G^n$ is isomorphic
to an extension of a finite subgroup $K$ of $G$ by a subgroup of
the direct product of two groups one of which is a dihedral group
of order $2n_1$ and the other is a direct product of several
copies of a dihedral group of order $2n_2$.
%\item[(d)] The subgroup $G^n$ is torsion free and $\bigcap_{k=1}^\infty G^{kn}=\{1\}$.
%\footnote{vsio nuzhno? pro punkty b,c,d}
\end{itemize}
\end{theorem}

\begin{rem} \label{rem:n}
As it was shown in \cite[Sec. 3]{Iv-Olsh}, %\footnote{tochnee ssylku}
the exponent $n$ from the previous statement can be chosen to be
any sufficiently large integer which is divisible by
$2^{k_0+5}n_0$, where
\begin{itemize}
\item[a)] $n_0/2$ is the least common multiple of the exponents of
the holomorphs $Hol(K)$ over all finite subgroups $K$ of $G$;
\item[b)]$k_0$ is the minimal integer with $2^{k_0-3} > \max |K|$
over all finite subgroups $K$ of $G$.
\end{itemize}
\end{rem}

The question which now naturally arises concerns the number (up to
isomorphism) of periodic quotients for a fixed $G$.

In \cite{Olsh-bulg} the second author of the present article
proved that there exists a continuum of $2$-generated
non-isomorphic groups with the property that every proper subgroup
is cyclic of prime order. In \cite{Der} G. Deryabina showed that
for every odd prime $p$ there is continuum of non-isomorphic
simple groups, generated by two elements, such that all their
proper subgroups are cyclic $p$-groups.

First examples of infinite finitely generated periodic residually finite groups were found by
E. Golod in \cite{Golod}. The techniques that he used allow, in fact, to obtain
uncountably many of such groups (we thank L. Bartholdi for pointing this out). Explicitly, constructions of
continuous families of finitely generated periodic residually finite groups first appear in the papers
\cite{Grigorchuk-84},\cite{Grigorchuk_p-gps} of R. Grigorchuk (see also
\cite[2.10.5]{Nekrashevych}).

Speaking of periodic groups of bounded exponent, one should
mention a result of V. Atabekyan \cite{Atab}, who demonstrated
that the free group $F_2$, of rank $2$, possesses a continuum of
pairwise non-isomorphic simple quotient-groups of given period $n$
for each sufficiently large odd $n$. The third author of the
present article established an analogous result in the situation
when $n$ is a sufficiently large even integer divisible by $2^9$ (see
\cite{Son}).
%(in fact, he proved a much stronger statement, claiming that for any such $n$ there is a
%continuum of non-isomorphic $2$-generated simple groups of exponent $n$ -- see \cite{Son}).

The first objective of this paper is to present several ways in
which one can construct periodic quotients of a given
non-elementary hyperbolic group $G$. We suggest three different
methods each of which gives a continuum of such non-isomorphic
quotients.
%\footnote{eshe fin.approksimiruemye primery?}

The second and the third approaches rely heavily on the techniques
and statements used by Ivanov and Olshanskii in their proof of
Theorem \ref{thm:Iv-Olsh} in \cite{Iv-Olsh}. The first approach
also uses their result but in a less thorough manner. These
approaches are used to produce the following three theorems.
(Recall that a group $H$ is said to be a {\it central extension}
of the group $K$ if there exists a subgroup $N \lhd H$, contained
in the center $Z(H)$ of $H$, such that $H/N \cong K$.)

\begin{thm} \label{thm:contin_2} Let $G$ be a non-elementary hyperbolic group and let an integer $n$
be chosen according to Remark \ref{rem:n}. Then $G$ possesses a
continuum of pairwise non-isomorphic periodic quotients $Q_{1i}$,
$i \in I_1$, of exponent $2n$, each of which is a central
extension of the group $G/G^n$.
\end{thm}

\begin{thm} \label{thm:contin_1} Let $G$ be a non-elementary hyperbolic group and let an integer
$n$ be chosen according to Remark \ref{rem:n}. Then $G$ possesses
a continuum of pairwise non-isomorphic periodic quotients
$Q_{2i}$, $i \in I_2$, of unbounded exponent, such that for every
$i \in I_2$ the order of an arbitrary element $g \in Q_{2i}$
divides $2^l n$ for some $l=l(g) \in \N$.
\end{thm}

Every non-elementary hyperbolic group $G$ contains a unique
maximal finite normal subgroup $E(G)$ (\cite[Prop. 1]{Olsh2}).
The quotient $\hat G=G/E(G)$, in addition to being
non-elementary and hyperbolic, also satisfies $E(\hat G)=\{1\}$.

\begin{thm} \label{thm:contin_3} %centerless
Assume that $G$ is a non-elementary hyperbolic group and denote $\hat
G=G/E(G)$. Choose a positive integer $\hat n$ arising after an
application of
Remark \ref{rem:n} to $\hat G$. %\footnote{ne ochen' blagozvuchno}.
Then $G$ has a continuum of pairwise non-isomorphic centerless
quotients of exponent $2\hat n$.
\end{thm}

Observe that if $Z(K)=\{1\}$ and $H_1$, $H_2$ are two central
extensions of $K$ then $H_1/Z(H_1) \cong K \cong H_2/Z(H_2)$.

We will show that if the original group $G$ is centerless, then so
is the group $G/G^n$. Thus the collections of quotients produced
by Theorem \ref{thm:contin_2} and Theorem \ref{thm:contin_3} are
significantly different because for any $i,j \in I_2$ one has
$Q_{1i}/Z(Q_{1i})\cong Q_{1j}/Z(Q_{1j})$.

% Def. of res. finite gps; Lackenby's result, Olshanskii-Osin proof
%The fourth , discussed in this article, starts with a large group $G$.
Our second objective is to obtain an analogue of the above results
for large groups. Recall that a group $G$ is said to be {\it
large} if it has a normal subgroup $N \lhd G$ of finite index, for
which there exists an epimorphism $\varphi:N \twoheadrightarrow
F$, where $F$ is a non-abelian free group. A recent result of M.
Lackenby \cite{Lackenby} asserts that if $G$ is a finitely
generated large group and $g_1,\dots,g_r \in G$, then the quotient
$G/\langle\langle g_1^n,\dots,g_r^n\rangle\rangle$ is large for
infinitely many $n \in \N$. Lackenby's proof is topological. A
different, combinatorial, proof of this result (not requiring the large group to
be finitely generated) was found by A. Olshanskii
and D. Osin in \cite{Olsh-Osin}; we will employ their methods to
prove Theorem \ref{thm:quot-large} below.

%Let $H$ be a group and let $P \subset \N$ be an arbitrary set of primes. We will say that $H$ is a $\pi$-group if
%every element $h \in H$ has a finite order $o(h)$ and for each $h \in H$ every prime divisor of $o(h)$ is in $\pi$.

\begin{thm} \label{thm:quot-large} Suppose that $G$ is a finitely generated group having a normal subgroup $N$, of finite index,
that maps onto a non-abelian free group. Let $p$ be an arbitrary
prime number. Then $G$ possesses $2^{\aleph_0}$ of periodic
residually finite pairwise non-isomorphic quotients $H_{j}$, $j
\in J$, such that for every $j \in J$, the natural image of $N$ in
$H_j$ is a $p$-group.
\end{thm}

The output of the previous theorem is very different from the
continuous families given by
Theorems~\ref{thm:contin_2},\ref{thm:contin_1},\ref{thm:contin_3}
as it produces periodic quotients with the additional property of
residual finiteness. We recall that an infinite periodic finitely
generated group of bounded exponent cannot be residually finite
due to E. Zelmanov's positive solution of the restricted Burnside
problem \cite{Zelmanov1}-\cite{Zelmanov2}.
In the authors' opinion, the periodic quotients from the claim of Theorem \ref{thm:contin_1} are not very likely to be residually finite.
For instance, if $G$ is a non-elementary hyperbolic group possessing no proper finite index subgroups (it was shown I Kapovich and
D. Wise \cite{Kapovich-Wise} and, independently, by the second author \cite{Olsh-bass-lub}, that the existence of such a group is equivalent
to the existence of a non-residually finite hyperbolic group, which is a well-known open problem), then no non-trivial quotient of $G$ can be
residually finite. Thus an application of Theorem \ref{thm:contin_1} to $G$ will
produce a continuum of non-(residually finite) periodic groups.

Let $G$ be a finitely generated group and $p \in \N$ be a prime
number. Define a series of finite index characteristic subgroups
$\delta^p_i(G)$ of $G$ by $\delta^p_0(G)=G$ and
$$\delta^p_i(G)=\left[\delta^p_{i-1}(G),\delta^p_{i-1}(G)\right] \left(\delta^p_{i-1}(G)\right)^p
\le \delta^p_{i-1}(G),~i \in \N.$$ If $G$ is finitely generated,
then for every prime $p$ and $t\in \N$, the number
$|G/\delta^p_t(G)|$ is a natural invariant of $G$. Another point
in which the fourth approach differs from the others, is that it
uses these invariants in order to distinguish between the
quotients of $G$. Thus we show that for every prime $p$, there
exist continually many finitely generated residually finite
$p$-groups that can be distinguished via natural inner algebraic
invariants.

Let $\Omega=\{0,1\}^{\N}$ be the set of all infinite sequences of
$0$'s and $1$'s. In Section \ref{sec:large} we prove

\begin{cor} \label{cor:large_quot_explicit} Assume that $G$ is a finitely generated group, $N \lhd G$ is a
normal subgroup mapping onto $F_2$ and $G/N$ is a finite $p$-group
for some prime $p \in \N$. Then for every $\omega \in \Omega$, $G$
has quotient-group $H_\omega$ such that
\begin{itemize}
\item[(a)] $H_\omega$ is a periodic $p$-group; \item[(b)]
$\displaystyle \bigcap_{i=0}^\infty \delta^p_i (H_\omega)=\{1\}$,
thus $H_\omega$ is residually finite; \item[(c)] if $\omega\neq
\omega' \in \Omega$ then there is $t \in \N$ such that
$|H_\omega/\delta^p_t(H_\omega)|\neq
|H_{\omega'}/\delta^p_t(H_{\omega'})|$; consequently $H_\omega
\not\cong H_{\omega'}$.
\end{itemize}
\end{cor}

Finally, after combining Corollary \ref{cor:large_quot_explicit}
with a theorem of B. Baumslag and S. Pride~\cite{Baumslag-Pride},
we deduce
\begin{cor} \label{cor:more_gens_than_rels} Assume that a group $G$ admits a presentation with $k\ge 2$ generators and at most $k-2$ relators.
Then for any prime $p$, $G$ has a collection of quotients
$\{H_\omega~|~\omega \in \Omega\}$ such that
\begin{itemize}
\item $H_\omega$ is a periodic $p$-group for every $\omega \in
\Omega$; \item $H_\omega$ is residually finite for every $\omega
\in \Omega$; \item if $\omega\neq \omega' \in \Omega$ then there
is $t \in \N$ such that $|H_\omega/\delta^p_t(H_\omega)|\neq
|H_{\omega'}/\delta^p_t(H_{\omega'})|$; consequently $H_\omega
\not\cong H_{\omega'}$.
\end{itemize}
\end{cor}

%say something about explicit invariants in thm. 4?

%Note that since there are only countably many finitely generated
%groups with solvable word problem, most of the periodic groups
%given by Theorems
%\ref{thm:contin_2}-\ref{thm:quot-large} have
%unsolvable word and conjugacy problems. %\footnote{stoit li?}

%%%%%%%%%%%%%%%%%%%%%%%%%%%%%%%%%%%%%%%%%%%%%%%%%%%%%%%%%%%%%%%%%%%%%%%%%%%%%%%%%%%%%%
\section{A particular application}
%Explain how our four theorems can be applied to answer Wiegold's question.

Before proving Theorems \ref{thm:contin_2}-\ref{thm:quot-large} we
are going to consider one application. We will show that each of
these theorems can be used to answer the following question of J.
Wiegold :

\vspace{.2cm} {\bf Question.} (\cite[16.101]{Kourovka}) Does the
group $G=\langle x,y~\|~x^2=y^4=(xy)^8=1 \rangle$ have uncountably
many quotients each of which is a $2$-group?

\vspace{.2cm} As mentioned by Wiegold, the group $G$ can be
homomorphically mapped onto a $2$-group $\Gamma$ that is a
subgroup of index $2$ in the first Grigorchuk's group
$\mathcal{G}$ (introduced in \cite{Grig-first_gp}). It is known
that every proper quotient of $\mathcal{G}$ is finite
\cite{Grig-just_inf}, in other terms, $\mathcal G$ is {\it just infinite}.
In a private communication with the first
author, R. Grigorchuk asserted that a similar property holds for
$\Gamma$. Thus, the group $\Gamma$ has only countably many
distinct quotients. (More generally, L. Bartholdi has recently proved that every
finitely generated subgroup of Grigorchuk's group $\mathcal{G}$ has at
most countably many distinct quotients \cite{Bartholdi}.)

In order to apply our approaches to the group $G$, we need to
observe that $G$ is a subgroup of index $2$ in the group
$$T=\langle a,b,c~\|~a^2=b^2=c^2=1, (ab)^2=(bc)^4=(ac)^8=1\rangle.$$ And since $1/2+1/4+1/8<1$, $T$ is a reflection
group of isometries of the hyperbolic plane $\mathbb{H}^2$, whose
fundamental domain is a triangle with angles $\pi/2,\pi/4$ and
$\pi/8$. In fact, $G$ consists of all isometries from $T$
preserving the orientation on $\mathbb{H}^2$ (see
\cite[7.3.1]{Fine-Ros}).

Hence the group $G$ acts properly discontinuously and cocompactly
on $\mathbb{H}^2$, therefore it is non-elementary (word)
hyperbolic. It is well-known that every finite subgroup of $T$
lies inside of a finite parabolic subgroup \cite[Ex. 2d, p.
130]{Bourbaki}. Hence, it must be a subgroup of a dihedral group
$D(2k)$, of order $2k$, where $k=2,4$ or $8$. Since $G$ contains
no reflections, every finite subgroup of $G$ will be cyclic of
order dividing $8$. Now, if $K$ is a cyclic group of order $2^m$,
$m\in \N$, then $|Aut(K)|=2^{m-1}$, hence $Hol(K)=K\rtimes Aut(K)$
has order $2^{2m-1}$. By Remark~\ref{rem:n}, the exponent $n$ from
the assumptions of Theorems \ref{thm:contin_2} and
\ref{thm:contin_1} can be chosen as a sufficiently large power of
$2$. Therefore both of these theorems give an affirmative answer
to Wiegold's question.

By a classical theorem of E. Cartan \cite[p. 267]{Cartan} (see
also \cite[Ch. II.2, Cor.~2.8]{Bridson-Haefliger}) every finite
subgroup $K$ of isometries of $\mathbb{H}^2$ fixes a non-empty
subset that is convex in a strong sense: for any two distinct
points fixed by $K$, $K$ will fix the entire bi-infinite geodesic
in $\mathbb{H}^2$ passing through these points. The fixed-point
set of a normal subgroup is invariant under the action of the
entire group, hence the fixed-point set $\mbox{Fix}(E(G))$ of
$E(G)$ cannot consist of a single point (because $G$ is infinite
and acts properly discontinuously). By the same reason, it cannot
be a bi-infinite geodesic. Therefore, $E(G)$ fixes at least three
points in general position, hence $\mbox{Fix}(E(G))=\mathbb{H}^2$.
And since the action of $G$ is faithful, one can conclude that
$E(G)=\{1\}$.
%As $G \lhd T$ and $E(G)$ is characteristic in $G$, we have $E(G) \lhd T$. Hence $E(G) \le E(T)$,
%implying that $E(G)=\{1\}$.
Thus Theorem \ref{thm:contin_3} will produce a continuum of
pairwise non-isomorphic centerless quotients of $G$ of exponent
$2n$ (for the same $n$ as above).

It remains to show that Theorem \ref{thm:quot-large} can also be
applied to $G$. It is well-known that hyperbolic triangle groups,
such as $G$, are large (\cite[7.3.1]{Fine-Ros}). For our purposes,
however, we will need to exhibit a stronger property, namely, the
existence of a normal subgroup $N \lhd G$, mapping homomorphically
onto a non-abelian free group, such that $|G:N|=2^l$ for some $l
\in \N$.

Denote by $A=\langle x \rangle_2 * \langle y \rangle_4$ the free
product of cyclic groups of orders $2$ and $4$.
Let $B$ be the
cartesian subgroup of $A$, i.e., the kernel of the natural
homomorphism from $A$ onto the direct product $\langle x \rangle_2
\times \langle y \rangle_4$. The cartesian subgroup of a free product
is a free group (this follows from Kurosh subgroup theorem and the
fact that the cartesian subgroup trivially intersects the free
factors). Note that $B$ is freely generated by its elements
 $[x,y],[x,y^2]$
and $[x,y^3]$ (where $[u,v]=uvu^{-1}v^{-1}$). Set $C=B^2$ to be
the normal subgroup of $A$ generated by the squares of elements of
$B$. Then $|B:C|=8$, $(xy)^4=[x,y][x,y^2]^{-1}[x,y^3] \in B
\setminus C$ and $(xy)^8 \in C$. By Schreier's formula, $C$ is a
free group of rank $(3-1)8+1=17$. Since $|A:\langle xy \rangle
C|=8$, we can use \cite[Lemma 2.3]{Olsh-Osin} to show that
$D=\langle \langle (xy)^8 \rangle \rangle^A \lhd C$ is a normal
closure in $C$ of at most $8$ elements. Hence the subgroup $E=C/D
\lhd A/D=G$ has a presentation with $17$ generators and $8$
relators. Since $17-8 \ge 2$, by a result of Baumslag and
Pride~\cite{Baumslag-Pride}, for every sufficiently large $k\in
\N$, $E$ has a normal subgroup $N'$, of index $k$, such that $N'$
can be mapped onto a non-abelian free group. Of course, one can
take $k=2^{l'}$ for some $l' \in \N$.
%Thus, we found a subgroup $N' \le G$ of index $|G:N'|=2^{l'} \cdot |A:C|=2^{l'+6}$ such that $N' \twoheadrightarrow F$.
Observe that $N=\bigcap_{g \in G} gN'g^{-1} \lhd G$ has finite
index in $N'$, and, therefore, it can also be mapped onto
non-abelian free group. It remains to note that $E/N$ is a finite
$2$-group because $gN'g^{-1} \lhd E$ and $|E:gN'g^{-1}|=2^{l'}$,
for every $g \in G$. Consequently $|G:N|=|G:E|\cdot|E:N|=2^6 \cdot
|E:N|=2^l$ for some $l \in \N$.

Therefore an application of Theorem \ref{thm:quot-large} to $G$
and $N$ provides a continuum of pairwise non-isomorphic residually
finite quotients of $G$ each of which is a $2$-group.

One can now see that Wiegold's problem mentioned above can be
solved in a number different ways. In the course of writing this
paper, R. Grigorchuk informed the authors of yet another possible
method to solve this problem,
%that can be carried out
using the techniques developed in his paper \cite{Grigorchuk-84}.

%%%%%%%%%%%%%%%%%%%%%%%%%%%%%%%%%%%%%%%%%%%%%%%%%%%%%%%%%%%%%%%%%%%%%%%%%%%%%%%%%%%%%%

\section{Preliminaries}
Fix a group $G$ with a finite symmetrized generating set $\A$. If
$g \in G$, $|g|_\A$ will denote the length of a shortest word $W$
over $\A$ representing $g$ in $G$. This gives rise to the standard
left-invariant distance function $d(\cdot,\cdot)$ on $G$ defined
by $d(x,y) = |x^{-1}y|_\A$ for any $x,y \in G$. Afterwards this
can be extended to the metric $d(\cdot,\cdot)$ on the Cayley graph
{\ga} in the usual way. For subset a $Q$ of {\ga}, the closed
$\varepsilon$-neighborhood is defined by
$$\cN_\varepsilon(Q)=\{x\in \Gamma(G,\A)~|~\exists~ y \in Q \mbox{
such that } d(x,y) \le \varepsilon\}.$$ For any two points $x,y
\in \Gamma(G,\A)$, $[x,y]$ will denote a geodesic segment between
them.

Fix a number $\delta \ge 0$. A geodesic $n$-gon in {\ga} is called
$\delta$-slim if each of its sides is contained in the closed
$\delta$-neighborhood of the others. Using the definition of E.
Rips, we will say that $G$ is $\delta$-hyperbolic if every
geodesic triangle in its Cayley graph is $\delta$-slim. As a
consequence of this definition, it is easy to obtain that each
geodesic quadrilateral in {\ga} is $2\delta$-slim. Further on we
shall assume that the group $G$ is non-elementary and
$\delta$-hyperbolic for some given $\delta \ge 0$. There is a
number of other (equivalent up to changing $\delta$) definitions
of $\delta$-hyperbolicity -- see \cite{Mihalik}, for example.
Since some of the Lemmas quoted below utilize them, we will
suppose that our $\delta$ is sufficiently large so that $G$ also
satisfies all of these other definitions.

%def. of cyclically reduced words;...
The length of a word $W$ over the alphabet $\mathcal A$ will be
denoted by $\|W\|$. Suppose $W$ represents an element $g \in G$.
Then we define $|W|_\A=|g|_\A$. If $p$ is a path in the Cayley
graph, then $p_-$, $p_+$ and $\|p\|$ will denote its {starting
point}, {ending point} and the {length} respectively. In the case
when $p$ is a simplicial path, $lab(p)$ will stand for the word
written on $p$ and $p^{-1}$ -- for the { inverse} path to $p$,
that is, $p^{-1}_-=p_+$, $p^{-1}_+=p_-$ and
$lab(p^{-1})\equiv lab(p)^{-1}$ %\footnote{lab(p) est' tol'ko u celyh putey, sostoyaschih iz celyh reber grafa cayley}.
Given some numbers $\bar \lambda, \bar c$ satisfying $0< \lambda
\le 1$, $c \ge 0$, we will say that $p$ is $(\lambda, c)$-{\it
quasigeodesic} if for any subpath $q$ of $p$ one has $\lambda
\|q\| - c \le d(q_-,q_+)$. A word $W$ will be called $(\lambda,
c)$-{\it quasigeodesic} provided some (equivalently, any) path
$p$, with $lab(p)\equiv W$, is $(\lambda, c)$-{quasigeodesic} in
{\ga}.

The next two lemmas are standard properties of a hyperbolic metric
space:
\begin{lem} \label{lem:close} {\normalfont (\cite[5.6,5.11]{Ghys},\cite[3.3]{Mihalik})}
There is a constant $\nu=\nu(\delta,\lambda,c)$ such that for any
$(\lambda,c)$-quasigeodesic path $p$ in $\Gamma(G,{\mathcal A})$
and a geodesic $q$ with $p_- = q_-$, $p_+ = q_+$, one has $p
\subset {\cN}_\nu(q)$ and $q \subset {\cN}_\nu(p)$.
\end{lem}

\begin{lem}\label{lem:quadrangle}{\normalfont (\cite[Lemma 4.1]{paper2})} Consider a geodesic
quadrilateral $x_1x_2x_3x_4$ in the Cayley graph {\ga} with
$d(x_2,x_3)>d(x_1,x_2)+d(x_3,x_4)$. Then there are points $u,v \in
[x_2,x_3]$ such that $d(x_2,u) \le d(x_1,x_2)$, $d(v,x_3) \le
d(x_3,x_4)$ and the geodesic subsegment $[u,v]$ of $[x_2,x_3]$
lies $2\delta$-close to the side $[x_1,x_4]$.
\end{lem}

It is well known (see, for example, \cite{Olsh2}) that in a
hyperbolic group $G$ every element $g$ of infinite order is
contained in a unique maximal elementary subgroup $E_G(g)$, and
\begin{multline} \label{eq:E_G} E_G(g)=\{x \in G~|~xg^mx^{-1}=g^n~\mbox{for some}~ m,n \in \Z\setminus\{0\} \}=\\
\{x \in G~|~xg^nx^{-1}=g^{\pm n}~\mbox{for some}~ n \in \N \}.
\end{multline}
The subgroup $\displaystyle E_G^+(g)=\{x \in
G~|~xg^nx^{-1}=g^n~\mbox{for some}~ n \in \N \}$ has index at most
$2$ in $E_G(g)$.

%%%%%%%%%%%%%%%%%%%%%%%%%%%%%%%%%%%%%%%%%%%%%%%%%%%%%%%%%%%%%%%%%%%%%%%%%%%%%%%%%%%%%%%%%%%%%%%%%%%

\section{Central extensions}\label{sec:central}
Fix a presentation $\langle \A ~\|~R \in \R_0\rangle$ of the
initial non-elementary hyperbolic group $G=G(0)$, where the set of
generators $\A$ is finite and $\R_0$ is the set of all relators in $G$, that is, $\R_0$ consists of all words from
over the alphabet $\A^{\pm 1}$ that represent the identity element in $G$.

Choose the exponent $n$ by Theorem \ref{thm:Iv-Olsh} and consider
the group $G(\infty)=G/G^n$.

Let us recall a few things from \cite{Iv-Olsh}. It is shown that
\begin{equation} \label{eq:G_inf}
 G(\infty) = \langle \A ~\|~R, R \in \R_0, A_1^n,A_2^n,\dots\rangle,
\end{equation}
where the word $A_j$, $j \in \N$, is called the {\it period of
rank $j$}. Moreover, each group of rank $i$
\begin{equation*} \label{eq:G_i} G(i) \stackrel{def}{=} \langle \A ~\|~R, R\in\R_0, A_1^n,A_2^n, \dots, A_i^n\rangle
\end{equation*} is also non-elementary hyperbolic.
For every $i \in \N$, the word $A_i$ has infinite order in
$G(i-1)$ and there exists a unique maximal finite subgroup
$\F(A_i)$ which is normalized by $A_i$ in $G(i-1)$ (more
precisely, $\F(A_i)$ is the torsion subgroup of
$E^+_{G(i-1)}(A_i)$). One of the properties, proved in
\cite{Iv-Olsh}, states that $A_i^{n/2}TA_i^{-n/2}T^{-1}=1$ in
$G(i-1)$ for any $i \in \N$ and any word $T$ representing an
element of $\F(A_i)$ (see the adaptation of Lemma 18.5 from
\cite{Iv} in \cite[Section 14]{Iv-Olsh}).

A word $J$ is called an $\mathcal{F}(A_i)$-involution if $J$
normalizes the subgroup $\mathcal F(A_i)$ of $G(i-1)$, $J^2 \in
\mathcal F(A_i)$ and $J^{-1}A_iJ=A_i^{-1}T$ for some $T \in
\mathcal F(A_i)$. By Lemma 19.2 \cite{Iv-Olsh},\cite{Iv},
$A_i^{n/2}JA_i^{n/2}J^{-1}=1$ in the group $G(i-1)$. If for a
given $i \in \N$ there exists an $\mathcal{F}(A_i)$-involution,
the period $A_i$ is said to be {\it even}; otherwise, $A_i$ is
said to be {\it odd}. Equivalently, $A_i$ is odd if and only if
$E_{G(i-1)}(A_i)=E^+_{G(i-1)}(A_i)$.

Consider a diagram $\Delta$ over the presentation \eqref{eq:G_inf}
on some orientable surface $S$ (say, a plane or a sphere). Then
one can fix the same (say, clockwise) direction of the contours of
each cell in $\Delta$. If $\Pi$ is a cell corresponding to the
relation $A_j^n=1$, $j \in \N$, one defines its rank $r(\Pi)$ to
be $j$; by definition, $r(\Pi)=0$ if $\Pi$ corresponds to a
relator $R \in \R_0$. A vertex $o \in \partial\Pi$ is called a
{\it phase vertex} if, starting at $o$ and going along the contour
of $\Pi$ in the clockwise direction, we read the word $A_j^{\pm
n}$. The number $r(\Delta)=\max\{r(\Pi)~|~\Pi \in \Delta \}$ is
called the {\it strict rank} of the diagram $\Delta$. If the
strict rank of $\Delta$ is $i$, the {\it type} $t(\Delta)$ of
$\Delta$ is a sequence $(\tau_i,\tau_{i-1},\dots,\tau_0)$ where
$\tau_j$ is the number of cells of rank $j$ in $\Delta$. For the
convenience of inductive reasoning we will impose the {\it
short-lex} order on the set of all types.

Let $\Pi_k$, $k=1,2$, be two distinct cells of $\Delta$ of rank
$j$. The pair of cells $\Pi_1$ and $\Pi_2$ is said to be a {\it
reducible $j$-pair} if there is a simple path $t$ in the diagram
$\Delta$, connecting two phase vertices $o_1$ and $o_2$ on the
their boundary contours, such that one of the following holds:
\begin{itemize}
\item[1)] The period $A_j$ is even and the words written on
$\partial \Pi_1$ and $\partial \Pi_2$ starting with the vertices
$o_1$ and $o_2$ are the same (i.e., $A_j^{\pm n}$), and $lab(t)$
is an $\mathcal{F}(A_j)$-involution.

\item[2)] The words written on $\partial \Pi_1$ and $\partial
\Pi_2$ starting with $o_1$ and $o_2$ are mutually inverse (that
is, they are either $A_j^n$ and $A_j^{-n}$ or $A_j^{-n}$ and
$A_j^{n}$ respectively), and $lab(t)$ represents an element of
$\F(A_j)$ in $G(j-1)$.
\end{itemize}

%\begin{lem} \label{lem:even_prod} Suppose a word $W$ represents an element of $[F,N]N^2N_0$ in $F$. Then there exists a disk diagram
%$\Gamma$ over the presentation \eqref{eq:G_inf} whose boundary label (starting some vertex and reading in the clockwise direction)
%is letter-by-letter equal to $W$, and which has an even number of cells of each positive rank.
%\end{lem}
%
%\begin{proof} Clearly, by the definitions of the subgroups $[F,N]$ and $N^2$, one can write
%$W$ as a product $\prod_{t} U_tA_{j_t}^{\pm n}U_t^{-1} h$ where $h \in N_0$ and for every $l$ the sum of all exponents at
%$A_l$ is even. Therefore, as in the proof of van Kampen's Lemma (see, for example \cite{L-S}),
%one can construct a disk diagram $\Gamma$ over the presentation \eqref{eq:G_inf},
%whose boundary contour is labelled by $W$, which has an even number of cells of rank $l$ every $l \in N$.
%
%\end{proof}

In the case when the original group $G$ is free, the following
statement can be compared with \cite[Lemma 5.2]{Olsh-rost}.

\begin{lem} \label{lem:spher_diag} Let $(\tau_i,\tau_{i-1},\dots,\tau_0)$ be the type of
an arbitrary spherical diagram $\Delta$ over the presentation
\eqref{eq:G_inf}. Then for and any $k \in \{1,2,\dots,i\}$ the
integer $\tau_k$ is even.
%\footnote{dlya svobodnoi gruppy takoe utverzhdenie est' v \cite{Olsh-mono} (ne nashel...)}
\end{lem}

\begin{proof} Fix $k \in \{1,\dots,i\}$ and use induction on $t(\Delta)$. If $\Delta$ has no cells of positive rank the claim is trivial.
Otherwise, the adaptation of \cite[Lemma 6.2]{Iv} described in
Section 7 of \cite{Iv-Olsh} claims that $\Delta$ has a reducible
$j$-pair of cells $\Pi_1$ and $\Pi_2$ for some $1 \le j \le i$.
Let $t$ denote the corresponding simple path between some phase
vertices $o_1$ and $o_2$ of the contours $\partial \Pi_1$ and
$\partial \Pi_2$, and set $T\equiv lab(t)$. Denote by $\Gamma$ the
subdiagram of $\Delta$ consisting of $\Pi_1$, $\Pi_2$ and $t$.
Then the word written on the boundary of $\Gamma$ starting $o_1$
is $A_j^{\epsilon n} T A_j^{\hat \epsilon n} T^{-1}$ where
$\epsilon,\hat \epsilon \in \{1,-1\}$. By the definition of a
reducible $j$-pair and the properties, mentioned above, in
$G(j-1)$ we have $$A_j^{\epsilon n/2} T A_j^{\hat \epsilon n/2}
T^{-1}=1.$$

By van Kampen's Lemma there exists a disk diagram $\Gamma'$ over
the presentation \eqref{eq:G_inf} with $r(\Gamma') \le j-1$ whose
boundary label is $A_j^{\epsilon n/2} T A_j^{\hat \epsilon n/2}
T^{-1}$. Gluing $\Gamma'$ with a copy of itself along the two
different occurrences of $T$ on their boundaries one obtains a new
disk diagram $\Gamma''$ (with $r(\Gamma'')=r(\Gamma')$), whose
boundary label is letter-by-letter equal to $A_j^{\epsilon n} T
A_j^{\hat \epsilon n} T^{-1}$ which is the boundary label of
$\Gamma$. Obviously, the number $\chi_k$ of cells of rank $k$ in
$\Gamma''$ is twice that number for $\Gamma'$, and, hence, it is
even.

Now one can perform a standard diagram surgery, by cutting
$\Gamma$ out of $\Delta$ and replacing it with $\Gamma''$, to
obtain a new spherical diagram $\Delta'$. Since
$t(\Delta')<t(\Delta)$ we can use the induction hypothesis to show
that the number $\psi_k$ of cells of rank $k$ in $\Delta'$ is
even. Due to the construction, one has $\tau_k=\psi_k-\chi_k+2$
(if $k=j$) or $\tau_k=\psi_k-\chi_k$ (if $k \neq j$). In either
case $\tau_k$ will be even.
\end{proof}

Let $F$ be the free group on $\A$ and let $N_0$ and $N$ denote the
the kernels of the natural homomorphisms $F \to G$ and $F \to
G(\infty)$ respectively. The {\it mutual commutator} $[F,N]$ is
the normal subgroup of $F$ generated by all commutators of the
form $[x,y]$ where $x \in F$ and $y \in N$.
%As usual, $N^2$ will denote the subgroup of $F$ generated by the squares of elements of $N$.

\begin{lem} \label{lem:even_pwr} Suppose the word $W \equiv A_{i_1}^{n\tau_1 }A_{i_2}^{n\tau_2}\dots A_{i_s}^{n\tau_s}$,
$1\le i_1<i_2<\dots i_s$, represents an element of the subgroup
$[F,N]N^2N_0$ in $F$. Then the integers $\tau_1,\dots,\tau_s$ are
all even.
\end{lem}

\begin{proof} Clearly, by definitions of the subgroups $[F,N]$ and $N^2$, one can write
$W$ as a product $\prod_{t} U_tA_{j_t}^{\pm n}U_t^{-1} h$ where $h
\in N_0$, and for every $l$ the number of occurrences of $A_l^{\pm
n}$ in this product is even. Therefore, as in the proof of van
Kampen's Lemma (see, for example, \cite{L-S}), one can construct a
disk diagram $\Gamma$ over the presentation \eqref{eq:G_inf},
whose boundary contour is labelled by $W$ and which has an even
number $\psi_l$ of cells of every positive rank $l$. On the other
hand, the equality from the assumptions of the lemma gives rise to
a diagram $\Delta_2$ over the presentation \eqref{eq:G_inf} with
the same boundary contour and having exactly $\tau_k$ cells of
rank $i_k$ for each $k=1,2,\dots,s$. Gluing together $\Delta_1$
with a mirror copy of $\Delta_2$ along their boundary contours,
one will obtain a spherical diagram $\Delta_3$. Applying Lemma
\ref{lem:spher_diag} to $\Delta_3$ one gets that the number
$(\tau_k+\psi_{i_k})$ of the cells of rank $i_k$ in it is even.
Therefore $\tau_k$ must also be even.
\end{proof}

Set $L=G^n\lhd G$ and $M=[G,L]L^2\lhd G$; the full preimages of
these subgroups in $F$ are $N$ and $[F,N]N^2N_0$ respectively.
Then the group $G/M$ is a central extension of the group
$G(\infty)=G/L$ by the central subgroup $L/M$. The analogue of the
next theorem in the case when $G$ is a free group of finite rank
was proved in \cite[Thm.~5]{Olsh-rost}.

\begin{cor} \label{cor:centr_ext} The group $L/M$ is a direct product of the groups of order $2$ generated by the natural
images of $A_j^n$, $j\in \N$.
\end{cor}

\begin{proof} Evidently, $L/M$ is an abelian group generated by the images $\hat a_j$ of $A_j^n$, $j \in \N$,
and these images have order at most $2$ in it. Now, Lemma
\ref{lem:even_pwr} asserts that if an element of the form $\hat
a_{i_1}^{\tau_1} \cdots \hat a_{i_s}^{\tau_s}$ is trivial in
$L/M$, then $\tau_k$ is even for each $k=1,\dots,s$. Hence $L/M$
is a direct product of the subgroups $\langle \hat a_j\rangle$.
\end{proof}

The following observation is a consequence of the fact that there
can be at most countably many different homomorphisms from a given
finitely generated group $G$ to a fixed countable group.

\begin{rem} \label{rem:cont_quot_non-isom} Let $G$ be a finitely generated group, $I$ be a set of cardinality continuum and $\{N_i\}_{i\in I}$
be a family of normal subgroups of $G$ such that $N_i \neq N_j$
whenever $i \neq j$. Then the set of quotients $\{G/N_i\}_{i \in
I}$ contains continually many pairwise non-isomorphic groups.
\end{rem}

\begin{proof}[Proof of Theorem \ref{thm:contin_2}.] Clearly, the order of any element $g \in G/M$ divides $2n$.
It is known that an infinite hyperbolic group always contains an
element of infinite order (\cite[8.3.36]{Ghys}), hence $G(\infty)
\neq G(i)$ for every $i \in \N$. Consequently the set of periods
$\{A_i\}$ is infinite. Thus, according to Theorem
\ref{cor:centr_ext}, the abelian group $L/M$ has continually many
distinct subgroups. As $L/M$ is contained in the center of the
group $H=G/M$, $H$ also has a continuum of different central
normal subgroups $Z_i$, $i \in I$. Note that since $G$ is finitely
generated, then so is $H$. By Remark~\ref{rem:cont_quot_non-isom},
there is a continuum of non-isomorphic groups among the quotients
$\{H/Z_i~|~i\in I\}$, all of which are central extensions of
$G(\infty)$. Therefore, the statement of the theorem is true.
%Hence for every $i \in I$ the quotient
%$H/Z_i$ can be isomorphic to at most countably many other
%quotients \footnote{raz'yasnit'?} $H/Z_j$, $j\neq i$. Therefore
\end{proof}

%%%%%%%%%%%%%%%%%%%%%%%%%%%%%%%%%%%%%%%%%%%%%%%%%%%%%%%%%%%%%%%%%%%%%%%%%%%%%%%%%%%%%%%%%%%%%%%%%%%

\section{Quotients of unbounded exponent}\label{section:unbounded}

In the present section we will prove Theorem \ref{thm:contin_1}.
Suppose that the non-elementary hyperbolic group $G$ has a presentation $\langle \A ~\|~R \in \R_0\rangle$, where
$\A$ is finite and $\R_0$ is the set of all relators in $G$.

By \cite[Lemma 18]{Iv-Olsh}, there is a constant $k_0 \in \N$ such
that $2^{k_0-3} > \max|K|$ for any finite subgroup $K$ of $G$. By
$n_0$ denote a constant such that $n_0/2$ is the least common
multiple of the exponents of holomorphs $Hol(K)$ over all finite
subgroups $K$ of $G$. %Notice that $n_0$ is a power of $2$.

Let us recall the construction of the quotient $G/G^n$ from
\cite{Iv-Olsh},
% of a non-elementary hyperbolic group $G$ over the
%subgroup $G^n$ generated by $n$-th powers of elements of $G$,
 where $n$ is a large integer divisible by $2^{k_0+5}n_0$.

Introduce a total order $\prec$ on the set of words in the
alphabet $\A$ so that $\|X\|<\|Y\|$ implies $X \prec Y$. Set
$G(0)=G$ and define groups $G(i)$ by induction on $i$.

As a matter of notational convenience, in what follows we shall
often identify a word $W$ over the alphabet $\mathcal A$ with an
element that it represents in $G(0)$. For a word $W$, representing
an element of infinite order $w\in G(0)$, denote by $F(W)$ the
maximal finite subgroup of $E^+(W)$. Such a $W$ is called {\it
simple} in $G(0)$ (see \cite[Sec. 4]{Iv-Olsh}) if its coset in
$E_G(W)/F(W)$ generates the cyclic subgroup $E_G^+(W)/F(W)$ and
none of the conjugates of elements $W^{\pm 1}F$ ( where the word
$F$ represents an element from $F(W)$) has length less than
$\|W\|$.

Assuming that the group $G(i-1),~i \ge 1$ is already defined,
consider the least (with respect to the order $\prec$) word $A$ of
infinite order in $G(i-1)$, with additional requirement that $A$
is simple in $G(0)$ if $\|A\|<C$ (see \cite[Lemma 13]{Iv-Olsh} for
the meaning of parameter $C$). Declare such a word $A$ to be the
period $A_i$ of rank $i$, and define the group $G(i)$ by imposing
the relation $A_i^n=1$ on the group $G(i-1)$:
$$
G(i)=\langle\A~\|~R \in \mathcal R_0 \cup \{A_1^n,A_2^n, \dots,
A_i^n\}\rangle.
$$
The period $A_i$ exists for every $i \ge 1$, and the group $G/G^n$
is the direct limit of groups $G(i)$:
$$
G/G^n=G(\infty)=\langle\A~\|~R \in \mathcal R_0 \cup
\{A_1^n,A_2^n, \dots ~\}\rangle.
$$

The construction of Ivanov and Olshanskii \cite{Iv-Olsh} may be
modified in the following way. Instead of imposing on each step a
relation $A^n=1$ for a fixed exponent $n$, we will introduce
periodic relations with large exponent that may vary from step to
step. In the case when $G(0)$ a non-abelian free group such
modification has already appeared in \cite{Olsh-rost}.

Given a sequence $\omega=(\omega_j)_{j=1}^{\infty}$ of $0$'s and
$1$'s, we modify the choice of defining relations for the groups
$G(i)$ from \cite{Iv-Olsh} as follows.

Set $n(0)=n$ and, for $j \ge 1$, define $n(j)=(1+\omega_j)n(j-1)$.
Inductively define groups $G_{\omega}(i)$: set $G_{\omega}(0)=G$
and assuming that the period $A_i$ of rank $i$ is chosen, define
$G_{\omega}(i)$ to be the quotient of $G_{\omega}(i-1)$ by the
normal closure of $A_i^{n(i)}$:
$$
G_{\omega}(i)=\langle\A~\|~R \in \mathcal R_0 \cup
\{A_1^{n(1)},A_2^{n(2)}, \dots, A_i^{n(i)}\}\rangle.
$$
The direct limit of groups $G_{\omega}(i)$ with respect to
canonical epimorphisms is denoted by $G_{\omega}(\infty)$.

Note that for any sequence $\omega$ of $0$'s and $1$'s the
sequence $(n(j))$ is non-decreasing and $n(j)\ge n$ for every $j
\ge 0$. We now refer to \cite{Olsh-rost}, and modify the arguments
from \cite{Iv-Olsh} in a way similar to the one given in
\cite{Olsh-rost}. Namely, the order of the period $A_i$ in
$G_{\omega}(\infty)$ is now $n(i)$ instead of $n$; in all the
estimates, the terms $\|A_i^n\|$ and $n\|A_i\|$ are substituted by
$\|A_i^{n(i)}\|$ and $n(i)\|A_i\|$ respectively; finite subgroups
of $G_{\omega}(i)$ are isomorphically embedded into a direct
product of a direct power of the dihedral group $D(2n(i))$
(instead of $D(2n)$) and an elementary associated with $G$ group
(see \cite{Iv-Olsh}); $n$ is replaced by $n(i)$ in the group
identities of the analogue of Lemma 15.10
\cite{Iv-Olsh},\cite{Iv}; in the inductive step from
$G_{\omega}(i)$ to $G_{\omega}(i+1)$ (in arguments from
 \cite[\S\S~18,19]{Iv} and their analogues from \cite{Iv-Olsh}) $n$ is
replaced by $n(i+1)$.

After all of these modifications, we obtain that for any sequence
$\omega$ the period $A_i$ exists for every $i \ge 1$ and the group
$G_{\omega}(\infty)$ is infinite and periodic.

Let $\omega$, $\omega'$ be two different sequences of $0$'s and
$1$'s. Then the kernels of the canonical homomorphisms $G \to
G_{\omega}(\infty)$ and $G \to G_{\omega'}(\infty)$ are different.
Indeed, let $j$ be the smallest index where $\omega$ and $\omega'$
differ, say $\omega_j=0$, $\omega'_j=1$. We remark here that for
all $i$, $0 \le i < j$, the sets of defining words of the groups
$G_{\omega}(i)$ and $G_{\omega'}(i)$ coincide in view of
minimality of $j$ and the fact that the order $\prec$ was fixed
a-priori. This means, in turn, that the same word $A_j$ is the
period of rank $j$ in both $G_{\omega}(\infty)$ and
$G_{\omega'}(\infty)$. By the analogue of Lemma 10.4
\cite{Iv-Olsh} for $G_{\omega}(\infty)$ (for
$G_{\omega'}(\infty)$), the order of the period $A_j$ in
$G_{\omega}(\infty)$ (in $G_{\omega'}(\infty)$) is equal to $n(j)$
($2n(j)$ respectively).

Thus, the groups $G_{\omega}(\infty)$ and $G_{\omega'}(\infty)$
are quotients of the finitely generated group $G$ by different
normal subgroups provided $\omega \ne \omega'$.
%and
%therefore the set of different presentations
%$\{G_{\omega}(\infty)\}_{\omega}$ is of cardinality continuum.
Remark \ref{rem:cont_quot_non-isom} permits us to conclude that
%Since there are only countably many different homomorphisms of the
%finitely generated group $F(\A)$ onto a fixed countable group, we
%conclude that
the set of pairwise non-isomorphic groups among
$\{G_{\omega}(\infty)\}_{\omega}$ is of cardinality continuum.

%%%%%%%%%%%%%%%%%%%%%%%%%%%%%%%%%%%%%%%%%%%%%%%%%%%%%%%%%%%%%%%%%%

\section{Aperiodic elements in hyperbolic groups}
The purpose of this section is to establish a few auxiliary
results, that will be used in Section~\ref{sec:hard} in order to
develop the third approach and prove Theorem~\ref{thm:contin_3}.
%The proofs of
%Theorems~\ref{thm:contin_2},\ref{thm:contin_1},\ref{thm:quot-large} do not use these technical statements
%and a reader interested only in the first, second or fourth approaches could just skip forward to the
%corresponding sections.

Throughout this section we will assume that $G$ is a
non-elementary hyperbolic group with a fixed finite symmetrized
generating set $\A$.

Let $W_1,W_2,\dots,W_l$ be words in $\A$ representing elements
$w_1,w_2,\dots,w_l$ of infinite order, where $E_G(w_i) \neq
E_G(w_j)$ for $i \neq j$. For any given $M\ge 0$ consider the set
$S(W_1,\dots,W_l;M)$ consisting of words
$$W \equiv W_{i_1}^{m_1}W_{i_2}^{m_2} \dots W_{i_s}^{m_s}$$ where
$s \in \N$, $i_k \neq i_{k+1}$ for $k=1,2,\dots,s-1$ (each $i_k$
belongs to $\{1,\dots,l\}$), and $|m_k|>M$ for $k=2,3,\dots,s-1$.
The following lemma will be useful:

\begin{lem}\label{lem:quasigeodesic}{\normalfont (\cite[Lemma 2.3]{Olsh2})}
There exist constants $\lambda_1 =
\lambda_1(W_1,W_2,\dots,W_l)>0$, $c_1=c_1(W_1,W_2,\dots,W_l) \ge
0$ and $M_1 = M_1(W_1,W_2,\dots,W_l)>0$ such that any path $p$ in
the Cayley graph {\ga} with $lab(p) \in S(W_1,\dots,W_l;M_1)$ is
$(\lambda_1,c_1)$-quasigeodesic.
\end{lem}

Consider a closed path $p_1q_1p_2q_2$ in the Cayley graph {\ga}
such that $lab(q_1)$, $lab(q_2^{-1}) \in S(W_1,\dots,W_l;M_2)$.
Thus $q_1=o_1\dots o_s$ %\footnote{bukva o ne nuzhna dlya vershin?}
where $lab(o_k)\equiv W_{i_k}^{m_k}$, $i_k \in \{1,\dots,l\}$,
$k=1,\dots,s$, and $i_k \neq i_{k+1}$, $k=1,\dots, s-1$.
Similarly, $q_2^{-1}=\bar o_1\dots \bar o_{\bar s}$ where
$lab(\bar o_k)\equiv W_{j_k}^{\bar m_k}$, $j_k \in \{1,\dots,l\}$,
$k=1,\dots,\bar s$, and $j_k \neq j_{k+1}$, $k=1,\dots, \bar s-1$.
We will say that $v$ is a phase vertex of a path $o_j$ if the
subpath of $o_j$ from $(o_j)_-$ to $v$ is labelled by some power
of the word $W_{i_j}$ (and similarly for $\bar o_j$).

Paths $o_k$ and $\bar o_{\bar k}$ will be called {\it compatible}
if there is a path $u_k$ in {\ga} joining some phase vertices of
$o_k$ and $\bar o_{\bar k}$ such that $lab(u_k)W_{{j_{\bar k}}}
lab(u_k)^{-1}$ represents an element of $E_G(w_{i_k})$ in $G$.
Such a path $u_k$ is said to be {\it matching}.

The following is a simplification of \cite[Lemma 2.5]{Olsh2} (we
note that instead of requiring $s$ and $\bar s$ to be bounded in
the proof of \cite{Olsh2} it is enough to demand that one uses
only finitely many distinct words $W_1,\dots,W_l$):

\begin{lem} \label{lem:lemma2.5-olsh2} Suppose that $q_1,q_2,p_1,p_2$ are as above, and
$\|p_1\|,\|p_2\| \le C_1$ for some $C_1$. Then there exist
integers $M_2$ and $\epsilon\in \{-1,0,1\}$, such that the paths
$o_k$ and $\bar o_{k+\epsilon}$ are compatible for all
$k=2,\dots,s-1$, whenever $|m_1|,\dots,|m_{s}| \ge M_2$, $|\bar
m_1|,\dots,|\bar m_{\bar s}| \ge M_2$.
\end{lem}

Recall (see \cite{Olsh2}) that two elements $g,h$ having infinite
order in $G$ are
said %\footnote{nazvany v [??]}
to be {\it commensurable} if there exist $a \in G$ and $k, l \in
\Z\setminus \{0\}$ such that $g^k=ah^la^{-1}$.
%In this case we will write $g \sim h$. Obviously
%"$\sim$" is an equivalence relation on the set of elements of infinite order in the group. The notation $g \nsim h$ will be used
%when $g$ and $h$ are not commensurable.
%
%\begin{lem} {\rm (\cite[Lemma 3.1]{Olsh2})} \label{lem:lemma3.1-olsh2} If elements $x$ and $y$ of a hyperbolic group
%have infinite order and satisfy $E_G(x)\neq E_G(y)$, then for any sufficiently large $|u|$ and $|v|$ the element $z=x^uy^v$
%has infinite order and is not commensurable with either $x$ or $y$.
%\end{lem}

Suppose, now, that the maximal finite normal subgroup $E(G)$ of a
non-elemen\-tary hyperbolic group $G$ is trivial. An element $g\in
G$ is said to be {\it suitable} if it has infinite order and
$E_G(g)=\langle g \rangle$.

\begin{lem} {\rm (\cite[Lemma 3.8]{Olsh2})} \label{lem:lemma3.8-olsh2} Every non-elementary hyperbolic
group $G$, such that $E(G)=\{1\}$, contains infinitely many
pairwise non-commensurable suitable elements.
\end{lem}

%\begin{lem} \label{lem:y-mod} {\rm (\cite[Lemma 4.3]{paper3})} Let $G$ be a non-elementary hyperbolic
%group with $E(G)=\{1\}$ and let $g \in G$ be a suitable element. If $y \in G \backslash E_G(g)$ then there exists $N \in \N$
%such that the element $yg^n$ is suitable for every $n\ge N$.
%\end{lem}
%
A subgroup $H \le G$ is called {\it malnormal} if $gHg^{-1} \cap
H=\{1\}$ for all $g\in G\setminus H$. In the special case when $G$
is torsion-free, the following lemma was proved by I.~Kapovich in
\cite[Thm. C]{Kapovich}.

\begin{lem} \label{lem:maln} Let $G$ be a non-elementary hyperbolic group with $E(G)=\{1\}$. Then there are words
$W_1$ and $W_2$, representing elements of infinite order $w_1$ and
$w_2$ in $G$, and constants $0<\lambda_1\le 1$, $c_1\ge 0$ such
that
\begin{itemize}
\item any path in {\ga} labelled by a word from $S(W_1,W_2;0)$ is
$(\lambda_1,c_1)$-quasi\-geodesic;
%of the form $ W_{i_1}^{m_1}W_{i_2}^{m_2} \dots W_{i_{s-1}}^{m_{s-1}}W_{i_s}^{m_s}$ where $i_k\meq i_{k+1}$,
\item the subgroup $H=\langle w_1,w_2 \rangle \le G$ is free of
rank $2$; \item $H$ is malnormal in $G$.
\end{itemize}
\end{lem}

\begin{proof} By Lemma \ref{lem:lemma3.8-olsh2} one can find four pairwise non-commensurable suitable elements $w,x,y,z \in G$
represented by some words $W,X,Y,Z$ over the alphabet $\A$
respectively. Apply Lemma \ref{lem:quasigeodesic} to the set of
words $S(W,X,Y,Z;M_1)$ to find the constants $\lambda_1$, $c_1$
and $M_1$. Find $\nu_1=\nu_1(\delta,\lambda_1,c_1)$ according to
Lemma \ref{lem:close} and denote $C_1=2\delta+2\nu_1$. Now, let
$M_2$ be the constant from the claim of Lemma
\ref{lem:lemma2.5-olsh2} applied to any closed path $p_1q_1p_2q_2$
where $lab(q_i) \in S(W,X,Y,Z;M_2)$ and $\|p_i\| \le C_1$,
$i=1,2$. Fix an arbitrary integer $m>\max\{M_1,M_2\}$ satisfying
$\lambda_1 m-c_1>0$, and define $w_1=w^mx^m$, $w_2=y^mz^m$. Note
that since the elements $w_1$ and $w_2$ are represented by the
words $W_1\equiv W^{m}X^m$ and $W_2 \equiv Y^{m}Z^m$ and $m>M_1$,
every word from $S(W_1,W_2;0)$ belongs to $S(W,X,Y,Z;M_1)$.
Therefore the first claim of the lemma holds.

To show that $w_1$ and $w_2$ freely generate the subgroup
$H=\langle w_1,w_2 \rangle \le G$ it is enough to check that any
non-empty word $A \in S(W_1,W_2;0)$ is non-trivial in $G$. But
this follows immediately from the
$(\lambda_1,c_1)$-quasigeodesity:
$$|A|_G\ge \lambda_1 \|A\| -c_1 \ge \lambda_1 \min\{ \|W_1\|,\|W_2\| \}-c_1\ge \lambda_1 m
-c_1>0.$$ %\footnote{esli eto rassuzhdenie ne ochen' nuzhno, to mozhno
%prosto dat' ssylku: esli w,x,... v poparno raznyh E(w), E(x), .. ,
%to $w^m$, $x^m$, .. svobodno porozhdayut svobodnuyu podgruppu pri
%dostaochno bol'shom m}

Arguing by contradiction, assume that $H$ is not malnormal, i.e.,
$ba_1b^{-1}=a_2$ for some $b \in G\setminus H$, $a_1,a_2 \in
H\setminus \{1\}$; then $ba_1^lb^{-1}=a_2^l$ for all $l \in \N$.
Let $A_1,A_2 \in S(W_1,W_2;0)\subset S(W,X,Y,Z;M_2)$ and $B$ be
some words representing $a_1,a_2$ and $b$. Observe that one can
assume that the words $A_1$ and $A_2$ are cyclically reduced (with
respect to $W_1$ and $W_2$) after replacing $B$ with some word
$B_1$ representing an element $b_1$ of the double coset $HbH$; one
will still have $b_1 \notin H$ because $HbH \cap H=\emptyset$.
Then the words $A_1^l$ and $A_2^l$ will also belong to
$S(W,X,Y,Z;M_2)$. Consider a closed path $p_1'q_1'p_2'q_2'$ in
{\ga} such that $lab(p_1')\equiv B_1$, $lab(q_1')\equiv A_1^l$,
$lab(p_2')\equiv B_1^{-1}$, $lab(q_2')\equiv A_2^{-l}$.

As the path $q_1'$ is $(\lambda_1,c_1)$-quasigeodesic, one can
take $l$ to be so large that it will have a subpath $q_1$ whose
endpoints will be at distances greater than $(\|B_1\|+\nu_1)$ from
the endpoints of $q_1'$ and which will be labelled by one of the
following words: $X^mY^m$ (coming from an occurrence of the
subword $W_1W_2$ in the middle of $A_1^l$), $X^mZ^{-m}$ (coming
from an occurrence of the subword $W_1W_2^{-1}$ in the middle of
$A_1^l$), $W^{-m}Y^m$ (from $W_1^{-1}W_2$), $W^{-m}Z^{-m}$ (from
$W_1^{-1}W_2^{-1}$), $X^mW^m$ (from $W_1W_1$), $W^{-m}X^{-m}$
(from $W_1^{-1}W_1^{-1}$), $Y^mZ^m$ (from $W_2W_2$) and
$Z^{-m}Y^{-m}$ (from $W_2^{-1}W_2^{-1}$). We shall consider only
the first case, when $lab(q_1)\equiv X^mY^m$, because the other
cases are completely similar. Thus, $q_1=o_1o_2$ where
$lab(o_1)\equiv X^m$, $lab(o_2) \equiv Y^m$.

By Lemma \ref{lem:close} there are points $\alpha_1,\alpha_2 \in
[(q_1')_-,(q_1')_+]$ such that $d(\alpha_1,(q_1)_-)\le \nu_1$ and
$d(\alpha_2,(q_1)_+) \le \nu_1$. Then by Lemma
\ref{lem:quadrangle}, applied to the geodesic quadrilateral with
vertices $(p_1')_-$, $(q_1')_-$, $(p_2')_-$ and $(q_2')_-$, one
can find points $\beta_1,\beta_2 \in [(p_1')_-,(q_2')_-]$
satisfying $d(\alpha_i,\beta_i)\le 2\delta$, %\footnote{pochemu?? pro dva delta}
$i=1,2$. Applying Lemma \ref{lem:close} once again one will obtain
points $\gamma_1,\gamma_2 \in q_2'$ with $d(\gamma_i,\beta_i) \le
\nu_1$. Let $q_2$ denote the subpath of $q_2'$ (or $q_2'^{-1}$)
starting with $\gamma_2$ and ending with $\gamma_1$, and
$p_1=[(q_2)_+,(q_1)_-]$, $p_2=[(q_1)_+,(q_2)_-]$. According to the
triangle inequality one has
$$\|p_1\|=d((q_2)_+,(q_1)_-) \le d((q_1)_-,\alpha_1)+d(\alpha_1,\beta_1)+d(\beta_1,\gamma_1)\le 2\delta+2\nu_1=C_1,$$
and, similarly, $\|p_2\|\le C_1$. Note that the labels of $q_1$
and $q_2^{-1}$ belong to the set $S(W,X,Y,Z;M_2)$, hence by Lemma
\ref{lem:lemma2.5-olsh2}, $q_2^{-1}$ contains subpaths $\bar o_1$
and $\bar o_2$ with $(\bar o_1)_+=(\bar o_2)_-$,
 $lab(\bar o_i) \in \{W^m,X^m,Y^m,Z^m\}^{\pm 1}$, $i=1,2$,
such that there exist matching paths $u_i$ which connect some
phase vertices of $o_i$ and $\bar o_i$, $i=1,2$. Since the
elements $w,x,y,z$ are pairwise non-commensurable, by the
definition of compatible paths one immediately gets $lab(\bar
o_1)\equiv X^{\epsilon_1 m}$ and $lab( o_2) \equiv Y^{\epsilon_2
m}$ for some $\epsilon_i \in \{-1,1\}$, $i=1,2$. Moreover
$\epsilon_i=1$, $i=1,2$, because $E_G(x)=E_G^+(x)=\langle x
\rangle$ and $E_G(y)=E_G^+(y)=\langle y \rangle$. Thus the subpath
$\bar o_1 \bar o_2$ of $q_2^{-1}$ (and, hence, of $q_2'^{\pm 1}$)
originates from an occurrence of the same subword $W_1W_2$ in
$lab({q_2'}^{\pm 1})$.

Now, observe that since the word $lab(u_1)X lab(u_1)^{-1}$
represents an element of $\langle x \rangle$ then, applying
\eqref{eq:E_G}, one achieves $lab(u_1) = X^{t_1}$ in $G$ for some
$t_1 \in Z$. Similarly, $lab(u_2) = Y^{t_2}$ in $G$ for some $t_2
\in \Z$. Since the endpoints of $u_1$ and $u_2$ are phase vertices
of $o_1$, $\bar o_1$ and $o_2$, $\bar o_2$, there are paths $r_i$,
$i=1,2$, having the same starting point $(r_i)_-=(o_1)_+=(o_2)_-$
and the same ending point $(r_i)_+=(\bar o_1)_+=(\bar o_2)_-$,
$i=1,2$, such that $lab(r_1)$ represents an element of $\langle x
\rangle$ and $lab(r_2)$ represents an element of $\langle y
\rangle$. But $\langle x \rangle \cap \langle y \rangle = \{1\}$
and $lab(r_1)=lab(r_2)$ in $G$, therefore $lab(r_1)$ represents
the identity element in $G$.

Let $q_1''$, $q_2''$ denote the subpaths of $q_1'^{-1}$ and
$q_2'^{-1}$ respectively, satisfying $(q_1'')_-=(r_1)_-$,
$(q_1'')_+=(p_1')_+$, $(q_2'')_-=(p_1')_-$ and
$(q_2'')_+=(r_1)_+$. By construction of the paths $o_1,o_2,\bar
o_1, \bar o_2$, the words $lab(q_i'')$, $i=1,2$, belong to
$S(W_1,W_2;0)$, and, hence, the elements which they represent in
$G$ belong to $H$. Finally, the equality $b_1=lab(p_1')=lab(q_2'')
lab(r_1^{-1}) lab(q_1'')$ combined with $lab(r_1^{-1})=1$ imply
that $b_1 \in H$, contradicting the initial assumption. Thus the
lemma is proved.
\end{proof}

A word $B$ over $\A$ (and the element $b \in G$ represented by it)
is said to be {\it cyclically reduced} if for any word $U$,
conjugate to $B$ in $G$, one has $\|B\|\le \|U\|$. For $m\in \Z$,
any subword of $B^m$ is called $B$-{\it periodic}.

The following fact was established in \cite[Lemma 27]{Olsh1} for
the case of a torsion-free hyperbolic group $G$ by Olshanskii; in
\cite[Lemma 12]{Iv-Olsh} Ivanov and Olshanskii noted that the
claim continues to hold even if $G$ possesses elements of finite
order.

\begin{lem}\label{lem:cyc_min-qd} Assume that $G$ is a non-elementary hyperbolic group with a finite generating set
$\mathcal{A}$. There exist numbers $0<\lambda_0\le 1$ and $c_0 \ge
0$ such that for every cyclically reduced word $B$ (over
$\mathcal{A}$), representing an element of infinite order in $G$,
any $B$-periodic word $V$ is
$(\lambda_0,c_0)$-quasigeodesic. %\footnote{G-neelementarnaya gyperbolicheskaya}
\end{lem}

\begin{lem} \label{lem:comm_short} Let $G$ be a $\delta$-hyperbolic group and $\bar \lambda, \bar c$ be some numbers satisfying
$0<\bar \lambda \le 1$ and $\bar c \ge 0$. Then there exists a
constant $\Lambda_1=\Lambda_1(\delta,\bar \lambda,\bar c)\ge 0$
such that the following holds.

Assume $X_1,X_2$ are some words over $\A$ having infinite order in
$G$ and for each $i=1,2$, any path in {\ga} labelled by a power of
$X_i$ is $(\bar \lambda, \bar c)$-quasigeodesic. If
$X_1^l=AX_2^kA^{-1}$ in $G$ for some word $A$ and non-zero
integers $k,l$, then one can find words $U,W$ and a $X_1$-periodic
word $V$ such that $X_2=UVW$ in $G$ and $\|U\|,\|W\|\le
\Lambda_1$.
\end{lem}

\begin{proof} Let
$\bar \nu=\bar \nu(\delta,\bar \lambda,\bar c)$ be obtained from
the claim of Lemma \ref{lem:close}. Set
$\Lambda_1=2\delta+2\bar\nu$ and denote by $a \in G$ the element
represented by the word $A$.

By the assumptions, $X_1^{ls}=AX_2^{ks}A^{-1}$ in $G$ for every $s
\in \Z$ (pick $s$ so that $ks>0$). Consider the quadrilateral with
sides $p_1, q_1, p_2$ and $q_2$ in {\ga} where
$(p_1)_{-}=(q_2)_+=1$, $lab(p_1)\equiv A$, $(q_1)_-=(p_1)_+$,
$lab(q_1) \equiv X_2^{ks}$, $(p_2)_{-}=(q_1)_+$, $lab(p_2)\equiv
A^{-1}$, $(q_2)_-=(p_2)_+$, $lab(q_2) \equiv X_1^{-ls}$. Since the
paths $q_1$ and $q_2$ are $(\bar \lambda, \bar c)$-quasigeodesic,
 if one takes $|s|$ to be sufficiently large (compared to $|a|_\A$), it
will be possible to find a subpath $r$ of $q_1$ labelled by $X_2$
whose endpoints are at distances at least $(|a|_\A+\bar \nu)$ from
the endpoints of $q_1$. According to Lemma \ref{lem:close} there
are points $u,v$ lying on the geodesic segment $[(q_1)_-,(q_1)_+]$
at distances at most $\bar \nu$ from $r_-$ and $r_+$ respectively.
Because of the choice of the subpath $r$ one can now apply Lemma
\ref{lem:quadrangle} to find points $u',v' \in [(q_2)_+,(q_2)_-]$
such that $d(u,u') \le 2\delta$ and $d(v,v')\le 2\delta$. Again by
Lemma \ref{lem:close}, there are points $u'',v'' \in q_2$ situated
at distances at most $\bar \nu$ from $u'$ and $v'$ respectively.
Hence $d(r_-,u''),d(r_+,v'')\le 2\bar \nu+2\delta=\Lambda_1$.

Now the claim of the lemma will hold if one denotes by $U$ and $W$
the words written on the geodesic paths $[r_-,u'']$ and
$[v'',r_+]$ respectively and by $V$ the word written on the
subpath of $q_2$ (or $q_2^{-1}$) starting with $u''$ and ending
with $v''$.
\end{proof}

%def. of non-commensurability; $(\Lambda,t)$-aperiodicity (don't forget to mention that the word $Z$ must have infinite order)
%Reformulate $(\Lambda,t)$-aperiodicity in geometric terms.
Suppose $\Lambda >0$ and $t\in \N$. Following \cite{Iv-Olsh} we
shall say that an element $g\in G$ is $(\Lambda,t)$-{\it periodic}
if there are words $U$, $V$ and $W$ over $\A$ such that $g = UVW$
in $G$, $\max\{|U|_\A, |V|_\A\} \le \Lambda$, and $V$ is a
$Z$-periodic word, for some cyclically reduced word $Z$ having
infinite order in $G$, with $\|V\| \ge t\|Z\|$. An element $h \in
G$ is called $(\Lambda, t)$-{\it aperiodic} provided for every
factorization $h=h_1h_2h_3$, where $h_1,h_2, h_3 \in G$, such that
$|h|_\A = |h_1|_\A + |h_2|_\A + |h_3|_\A$, the element $h_2$ is
not $(\Lambda, t)$-periodic. It is not difficult to see that, in
geometric terms, $g$ is $(\Lambda, t)$-aperiodic if and only if
for any geodesic path $p$ in {\ga}, whose label represents $g$ in
$G$, and for any path $q$ such that $q_-,q_+ \in \cN_\Lambda(p)$
and $lab(q)\equiv V$ is a $Z$-periodic word (for some cyclically
reduced word $Z$ having infinite order in $G$), one has
$\|q\|=\|V\|<t\|Z\|$.

Suppose that $W_1,W_2$ are some words over the alphabet $\A$
representing elements $w_1,w_2 \in G$. Let $F^+(W_1,W_2)$ denote
the free monoid generated by the words $W_1$ and $W_2$ (i.e., the
set of all positive words in $W_1$ and $W_2$). Assume that
for some $\lambda_1>0$ and $c_1 \ge 0$, any word from
$F^+(W_1,W_2)$ is $(\lambda_1,c_1)$-quasigeodesic in $G$ and that
the canonical map $\psi:F^+(W_1,W_2) \to G$ is injective. Let $F^+
\subset G$ denote the set of elements represented by the words
from $F^+(W_1,W_2)$, i.e., $F^+=\psi(F^+(W_1,W_2))$. Note that
each non-trivial $f\in F^+$ has infinite order in $G$.

%Since the group $G$ is non-elementary, we can pick two elements $x,y\in G$ of infinite order with $E_G(x) \neq E_G(y)$.
%Choose some words $X,Y$ over $\A$ representing $x,y$.
%Upon replacing $X$ and $Y$ with their sufficiently large powers, we can make
%sure that any word from the free
%monoid $F^+(X,Y)$ is $(\lambda_1,c_1)$-quasigeodesic in $G$ where $\lambda_1$ and $c_1$ depend only on
%these words (see \cite[Lemma 2.3]{Olsh2}) and the canonical map $\psi:F^+(X,Y) \to G$ is injective (\cite[Cor. 6]{Olsh2}).
%Let $F^+ \subset G$ denote the set of elements represented by the
%words from $F^+(X,Y)$, i.e., $F^+=\psi(F^+(X,Y))$. Note that each non-trivial $f\in F^+$ has infinite order in $G$.
%

\begin{lem} \label{lem:aper_non-comm} %Suppose $F=\langle x,y \rangle\le G$ is a free quasiconvex subgroup
%of rank $2$ and let $F^+$ denote the set of positive elements (i.e., represented by positive words in
%$\{x,y\}$) of $F$.
There is a constant $\Lambda_1\ge 0$ such that for any
$\Lambda'\ge\Lambda_1$ there exist an integer $t'>0$ and
infinitely many pairwise non-commensurable elements in $F^+$ which
are $(\Lambda',t')$-aperiodic in $G$.
\end{lem}

\begin{proof} Let $\lambda_0$ and $c_0$ be the constants from Lemma \ref{lem:cyc_min-qd}.
Define $\bar \lambda = \min\{\lambda_0,\lambda_1\}$ and $\bar c =
\max\{c_0,c_1\}$, find the corresponding constant
$\Lambda_1=\Lambda_1(\delta,\bar \lambda,\bar c)\ge 0$ according
to Lemma \ref{lem:comm_short}, and take an arbitrary $\Lambda' \ge
\Lambda_1$.

In \cite[Lemmas 30,31]{Olsh1} it was shown that there exists $t'
\in \N$ and an infinite subset $\mathcal{B}=\{b_1,b_2,\dots\}
\subset F^+$ consisting of $(\Lambda',t')$-aperiodic elements in
$G$ (the proof works in the general case of an arbitrary
non-elementary hyperbolic group as observed in \cite[Lemma
15]{Iv-Olsh}).

For each $j \in \N$ find the (unique) word $B_j \in F^+(W_1,W_2)$
with $\psi(B_j)=b_j$. As $W_1$ and $W_2$ are words over $\A$, so
is $B_j$. Since $\mathcal{B}$ is infinite, after passing to its
infinite subset, we can assume that
$|b_{j+1}|_\A>t'\|B_j\|+2\Lambda_1$ for every $j\in \N$.

Suppose that there are indices $i<j$ such that $b_i$ is commensurable
with $b_j$. Choose a shortest representative $b$ of the conjugacy
class of $b_i$ in $G$; then $|b|_\A \le |b_i|_\A$. By the
assumption, there exists $a \in G$, $k \in \N$ and $l \in \Z
\backslash \{0\}$ such that $b^l=ab_j^ka^{-1}$.

Let $A$ and $B$ be shortest words over the alphabet $\A$
representing $a$ and $b$ respectively (note that the word $B$ is
cyclically reduced in $G$ by construction). Then
$B^{l}=AB_j^{k}A^{-1}$ in $G$ and one can use Lemma
\ref{lem:comm_short} to find words $U,V$ and $W$ such that
$B_j\stackrel{G}{=} UVW$ where $\|U\|,\|W\| \le \Lambda_1 \le
\Lambda'$ and $V$ is a $B$-periodic word. The latter equality
implies that
$$\|V\|\ge |V|_\A \ge |B_j|_\A-\|U\|-\|W\| \ge |b_j|_\A -2\Lambda_1> t'\|B_i\| \ge t'\|B\|,$$
which contradicts the $(\Lambda',t')$-aperiodicity of $b_j$.
Therefore no two distinct elements of $\mathcal B$ can be
commensurable and the statement is proved.
\end{proof}

\begin{lem} \label{lem:comm->aper} Suppose $G$ is a $\delta$-hyperbolic group,
$0<\bar \lambda \le 1$, $\bar c \ge 0$ and $\varkappa>0$. For any
$\Lambda>0$ there exists $\Lambda'>0$ such that for any $t' \in
\N$ there is $t\in \N$ satisfying the following.

Assume $X_1,X_2$ are some words over $\A$ representing elements of
infinite order $x_1,x_2 \in G$ and for each $i=1,2$, any path in
{\ga} labelled by a power of $X_i$ is $(\bar \lambda, \bar
c)$-quasigeodesic. If $x_1$ is $(\Lambda',t')$-aperiodic, $\|X_2\|
\le \varkappa \|X_1\|$ and $x_1^l=ax_2^ka^{-1}$ for some $a\in G$,
$k,l\in \Z\setminus \{0\}$, then $x_2$ is $(\Lambda,t)$-aperiodic.
\end{lem}

\begin{proof} Let $\lambda_0,c_0$ be from the claim of Lemma \ref{lem:cyc_min-qd},
$\nu_0=\nu_0(\delta,\lambda_0,c_0)$,
$\bar\nu=\bar\nu(\delta,\bar\lambda,\bar c)$ be as in Lemma
\ref{lem:close} and $\Lambda_1=2\delta+2\bar \nu$ be from Lemma
\ref{lem:comm_short}. Set $\Lambda'=4\delta+3\bar \nu +
\Lambda_1+\nu_0+\Lambda$, choose an arbitrary $t' \in \N$ and a
positive integer $m$ so that
\begin{equation*} \label{eq:m}
m>\frac1{\bar \lambda} (\varkappa+2\Lambda_1+\bar c) +2,
\end{equation*}
and denote $t=mt'$. Assume, now, that $x_1$ is
$(\Lambda',t')$-aperiodic.

Arguing by contradiction, suppose that $x_2$ is not
$(\Lambda,t)$-aperiodic. Then there is a geodesic path $p$ between
$1$ and $x_2$, points $a,b \in p$ and a
$(\lambda_0,c_0)$-quasigeodesic path $r$ in {\ga} such that the
label of $r$ is a $Z$-periodic word (for some cyclically reduced
word $Z$ over $\A$ having infinite order in $G$), $\|r\|\ge t
\|Z\|$ and $d(a,r_-),d(b,r_+) \le \Lambda$.

Using the assumptions together with Lemma \ref{lem:comm_short},
one can find a $(\bar \lambda,\bar c)$-quasi\-geodesic path $q$ in
{\ga} where $lab(q)$ is $X_1$-periodic and
$d(p_-,q_-),d(p_+,q_+)\le \Lambda_1$.

Since the geodesic quadrilaterals in {\ga} are $2\delta$-slim one
obtains
$$a,b \in p \subset \cN_{2\delta}([p_-,q_-]\cup [q_-,q_+] \cup [q_+,p_+]) \subset \cN_{2\delta+\Lambda_1}([q_-,q_+])
\subset \cN_{2\delta+\Lambda_1+\bar \nu}(q).$$

Thus, there are points $a',b' \in q$ with $d(a,a'),d(b,b') \le
2\delta+\Lambda_1+\bar \nu$; consequently $d(r_-,a'),d(r_+,b') \le
2\delta+\Lambda_1+\bar \nu+\Lambda$. Using $2\delta$-slimness of
the geodesic quadrilateral with the vertices $r_-,r_+,a',b'$
together with the property of quasigeodesics given by Lemma
\ref{lem:close}, one gets
\begin{multline}\label{eq:r_close} r \subset \cN_{\nu_0}([r_-,r_+]) \subset \cN_{\nu_0+2\delta}([r_-,a']\cup [a',b'] \cup [b',r_+])
\subset \\ \cN_{\nu_0+2\delta+2\delta+\Lambda_1+\bar
\nu+\Lambda}([a',b']) \subset \cN_{\nu_0+4\delta+\Lambda_1+\bar
\nu+\Lambda+\bar \nu}(q) = \cN_{\Lambda'-\bar\nu}(q).
\end{multline}

Now let us estimate the length of $q$ using its quasigeodesity and
the triangle inequality:
\begin{multline*} \|q\| \le \frac1{\bar \lambda} \bigl(d(q_-,q_+)+\bar c\bigr)\le
\frac1{\bar \lambda} \bigl(d(1,x_2)+d(1,q_-)+d(q_+,x_2)+\bar c\bigr) \le\\
\frac1{\bar \lambda} \bigl(\|X_2\|+2\Lambda_1+\bar c\bigr)\le
\frac1{\bar \lambda}\bigl(\varkappa \|X_1\|+2\Lambda_1+\bar
c\bigr).
\end{multline*}

Obviously, there is a path $\hat q$ in {\ga} such that $q$ is a
subpath of $\hat q$, $\|\hat q\|\le \|q\|+2\|X_1\|$ and $lab(\hat
q)\equiv X_1^s$ for some integer $s$. One has $$|s|\|X_1\| =
\|\hat q\| \le \|q\|+2\|X_1\|, ~\mbox{ hence } |s|\le \frac1{\bar
\lambda}\bigl(\varkappa +2\Lambda_1+\bar c\bigr)+2<m$$ due to the
choice of $m$. Since $t=mt'$, one can split $r$ in a concatenation
of $m$ subpaths $r=r_1 \dots r_m$ such that $lab(r_j)$ is a
$Z$-periodic word and $\|r_j\| \ge t' \|Z\|$ for all
$j=1,\dots,m$. The path $\hat q$, on the other hand, is a
concatenation of $|s|$ subpaths $q=\hat q_1\dots \hat q_{|s|}$,
each of which is labelled by $X_1^{\pm 1}$. Now, since $m>|s|$,
one can use the Pigeon-hole Principle to find $j,j' \in
\{1,\dots,m\}$, $j\le j'$, and $k \in \{1,\dots,|s|\}$, such that
$(r_j)_-,(r_{j'})_+ \in \cN_{\Lambda'-\bar \nu}(\hat q_k)\subset
\cN_{\Lambda'}([(\hat q_k)_-,(\hat q_k)_+])$ (here we utilize the
$(\Lambda'-\bar \nu)$-proximity of $r$ to $\hat q$ given by
\eqref{eq:r_close}). Let $r'$ be the subpath of $r$ starting at
$(r_j)_-$ and ending at $(r_{j'})_+$. Then $lab(r')$ is a
$Z$-periodic word and $\|r'\| \ge t' \|Z\|$. But this contradicts
the $(\Lambda',t')$-aperiodicity of $x_1$. The proof is finished.
\end{proof}

%def. of simple words in $G=G(0)$.

\begin{lem} \label{lem:aper_simple} Let $G$ be a non-elementary hyperbolic group with $E(G)=\{1\}$. For every $\Lambda>0$ there is
$t\in \N$ such that $G$ contains an infinite set $\{d_i~|~i \in
\N\}$ of pairwise non-commensurable $(\Lambda,t)$-aperiodic
elements of infinite order which are cyclically reduced and such
that $E_G(d_i)=\langle d_i \rangle$ for every $i \in \N$.
\end{lem}

\begin{proof} Let $W_1$, $W_2$, $H$, $\lambda_1$ and $c_1$ be as in Lemma \ref{lem:maln}, and
let $\lambda_0,c_0$ be from Lemma \ref{lem:cyc_min-qd}. Set $\bar
\lambda = \min\{\lambda_0,\lambda_1\}$, $\bar c =
\max\{c_0,c_1\}$, $\varkappa =1$ and find $\Lambda'>0$ from the
claim of Lemma \ref{lem:comm->aper}. Denote
$\Lambda''=\max\{\Lambda'+\nu_1,\Lambda_1\}$, where $\Lambda_1$
and $\nu_1=\nu_1(\delta,\lambda_1,c_1)$ are given by Lemmas
\ref{lem:aper_non-comm} and \ref{lem:close}.

According to Lemma \ref{lem:aper_non-comm}, there are $t'\in \N$
and an infinite set $\mathcal{B}=\{b_1,b_2,\dots\} \subset F^+
=\psi(F^+(W_1,W_2)) \subset H$ of pairwise non-commensurable
$(\Lambda'',t')$-aperiodic elements of infinite order in $G$.
 Observe that by \eqref{eq:E_G}, for any $h \in H \setminus \{1\}$
the malnormality of $H$ implies that $E_G(h)\subset H$, hence
$E_G(h)$ is cyclic as any torsion-free elementary group. In
particular, for each $j \in \N$ there is $b'_j \in F^+$ such that
$E_G(b_j)=\langle b'_j \rangle$ . For each $j\in \N$ one can
choose a word $B'_j \in F^+(W_1,W_2)$ representing $b'_j$ such
that any path in {\ga} labelled by a power of $B'_j$ is
$(\lambda_1,c_1)$-quasigeodesic (therefore it will also be $(\bar
\lambda, \bar c)$-quasigeodesic), because such a power will still
belong to $F^+(W_1,W_2)\subset S(W_1,W_2;0)$. Now, since
$b_j=(b_j')^s$ for some $s \in \N$, Lemma \ref{lem:close}
immediately implies that $b'_j$ is
$(\Lambda''-\nu_1,t')$-aperiodic because $b_j$ is
$(\Lambda'',t')$-aperiodic. Thus, the infinite set
$\mathcal{B}'=\{b'_1,b'_2,\dots\}$ consists of
$(\Lambda',t')$-aperiodic pairwise non-commensurable elements of
infinite order in $G$, with the additional property that
$E_G(b'_j)=\langle b'_j \rangle$.

For every $b'_j\in \mathcal{B}'$ choose a cyclically reduced
element $d_j \in G$ and an element $a_j \in G$ such that
$b'_j=a_jd_ja_j^{-1}$. Note that the set
$\mathcal{D}=\{d_1,d_2,\dots\}$ is an infinite set of pairwise
non-commensurable elements of infinite order in $G$, and
$E_G(d_j)=a_j^{-1}E_G(b'_j)a_j = \langle d_j \rangle$ is cyclic
for every $j \in \N$.

Let $t$ be given by Lemma \ref{lem:comm->aper}. Choose a shortest
word $D_j$ representing $d_j$ in $G$. Then $\|D_j\| \le \|B_j\|
=\varkappa \|B_j\|$ and, hence, according to Lemma
\ref{lem:cyc_min-qd} and the construction, one can apply Lemma
\ref{lem:comm->aper} to show that $d_j$ is $(\Lambda,t)$-aperiodic
for every $j \in \N$.
\end{proof}

%%%%%%%%%%%%%%%%%%%%%%%%%%%%%%%%%%%%%%%%%%%%%%%%%%%%%%%%%%%%%%%%%%%%%%%%%%%%%%%%%%%%%%%%%%%%%%%%%%%%%%%%%%%%%%%%%%%%%%%%%%%%%%%%%%%%%%%%

\section{Quotients of bounded exponent}\label{sec:hard}

In this section we are going to prove Theorem \ref{thm:contin_3}.
We will utilize the construction from \cite{Iv-Olsh} and explain
how to obtain series of infinite quotients of a non-elementary
hyperbolic group $G$ of bounded exponent by adding periodic
relations $A^{n(i)}$ with different exponents $n(i)$.

Let $G$ be a non-elementary hyperbolic group. Fix a presentation $\langle \A ~\|~R \in \R_0\rangle$ of $G$ as
in the beginning of Section \ref{section:unbounded}.
Since $E(G)$ is the maximal finite normal subgroup of $G$,
the quotient $\hat G=G/E(G)$ is again a non-elementary hyperbolic group with the additional property
$E(\hat G)=\{1\}$. After replacing $G$ with $\hat G$, we will further assume that $E(G)=\{1\}$.

Let $\nu_0=\nu_0(\lambda_0,c_0)$ be the constant provided by Lemma \ref{lem:close}, where
the pair $(\lambda_0,c_0)$ is chosen according to Lemma \ref{lem:cyc_min-qd}.
Choose the auxiliary parameters (such as $\Lambda, \alpha, \beta, \gamma,
\varepsilon, \rho, \theta, n$, etc.) as in \S~3 of \cite{Iv-Olsh}.
We refer to \cite{Iv-Olsh}, \cite{Iv} for the values and estimates
of the parameters; in the present section we will explicitly use
the following inequalities: $\theta>\alpha+3\beta$,
$2\gamma+1/2<\rho$.
%\footnote{ya dumayu, chto nado by otmetit', primernye znacheniya dlya $\alpha$,
%$\beta$, $\theta$ i $n$, a takzhe svyaz' mezhdu nimi. Dalee v
%dokazatel'stvah my pol'zuemsya neravenstvami $\alpha>1/2$,
%$\theta>\alpha+3\beta$, i t.p.}.
The choice of the parameter $n$, which is a very large integer, is made after all the other parameters are chosen.
Increasing $n$ (that is, multiplying it by a power of $2$),
if necessary, and using Lemma~\ref{lem:aper_simple}, %and \ref{lem:comm->aper}, %, \ref{lem:comm->aper} and 26 \cite{Olsh1}
we find an infinite set $\mathcal B=\{B_1,~B_2, \dots\}$ of
cyclically reduced words satisfying the following properties:

\begin{itemize}
\item[1.] Elements of the group $G$ represented by words
$B_1,~B_2, \dots $ are pairwise non-commensurable elements of
infinite order that generate their respective (cyclic) elementary subgroups
$E_G(B_k)$, $k=1,~2,~\dots$.

\item[2.] For every $B \in \mathcal B$, the word $B^{\pm 3}$
represents a $(\Lambda+\nu_0, \beta n)$-aperiodic element.

\end{itemize}

Passing to a subsequence, we can assume that $\|B_1\|>n^2$ and
$\|B_{k+1}\|>n^2\|B_k\|$ for $k \ge 1$. Given a sequence
$\omega=(\omega_k)_{k=1}^{\infty}$ of $0$'s and $1$'s, we are
going to construct a periodic quotient $G_\omega(\infty)$ of $G$
of exponent $2n$.

We set $G_{\omega}(0)=G$, introduce a total order $\prec$ and
define simple words as in Section~\ref{section:unbounded}.
Assuming that the presentation of $G_{\omega}(i-1)$ has already
been constructed, choose the period $A_i$ (of rank $i$) as in
Section~\ref{section:unbounded}.

We distinguish the following two cases:

\begin{itemize}
\item[(1)] A non-trivial power of some word $B_k \in \mathcal B$
is conjugate in $G_{\omega}(i-1)$ to an element $A_i^{\ell}F$, for
some $k \ge 1, \ell \in \mathbb Z, \ell \ne 0, F \in
\mathcal{F}(A_i)$. (It follows from Lemma \ref{lem:aper_non-comm},
Lemma \ref{lem:inforder_noncom} below and Lemma 18.5(c)
\cite{Iv},\cite{Iv-Olsh} in rank $i-1$ that such $k$ is unique.)
In this case we shall say that the period $A_i$ is {\it special}.

\item[(2)] None of non-trivial powers of words from $\mathcal B$
is conjugate in $G_{\omega}(i-1)$ to an element of the form
$A_i^{\ell}F$, for some $\ell \in \mathbb Z, \ell \ne 0, F \in
\mathcal F(A_i)$.

\end{itemize}

Set $n(i)=n$ if case (1) holds and $\omega_k=0$. Otherwise set
$n(i)=2n$.

The group $G_{\omega}(i)$ is obtained by imposing on
$G_{\omega}(i-1)$ the relation $A_i^{n(i)}=1$:
$$
G_{\omega}(i)=\langle\A~\|~R \in \mathcal R_0 \cup
\{A_1^{n(1)},A_2^{n(2)}, \dots, A_i^{n(i)}\}\rangle.
$$

Analysis of groups $G_{\omega}(i)$ is done using geometric
interpretation of deducing consequences from defining relations
according to the scheme of \cite{Iv-Olsh} (and \cite{Iv}). We
refer the reader to \cite{Iv-Olsh}, \cite{Iv} for the definitions
of a bond, a contiguity subdiagram and its standard contour, a
degree of contiguity, a reduced diagram, (strict) rank and type of
a diagram, a simple word in rank $i$ (some of these concepts were
defined in Section \ref{sec:central}).

References to lemmas from \cite{Iv-Olsh}, as well as lemmas from
\cite{Iv} and their analogs in \cite{Iv-Olsh}, will be made
preserving the numeration of \cite{Iv-Olsh} and of \cite{Iv}.

We note that our modification of the construction implies obvious
changes in formulations and proofs from \cite{Iv-Olsh},\cite{Iv}:
a cell of rank $i$ now corresponds to the relation $A_i^{n(i)}=1$
and $n(i)$ appears instead of $n$ in all the estimates of the
length of its boundary; the order of the period $A_i$ in
$G_{\omega}(\infty)$ is now $n(i)$ instead of $n$; finite
subgroups of $G_{\omega}(i)$ are isomorphically embedded into a
direct product of a direct power of the dihedral group $D(4n)$
(instead of $D(2n)$) and an elementary associated with $G$ group
(see \cite{Iv-Olsh}); in the inductive step from $G_{\omega}(i)$
to $G_{\omega}(i+1)$ (in arguments from \cite[\S\S~ 18,19]{Iv} and
their analogues from \cite{Iv-Olsh}) $n$ is replaced by $n(i+1)$.
The crucial part that needs explanation is validity of equations
in Lemma 15.10 \cite{Iv},\cite{Iv-Olsh}.

The following lemma about structure of diagrams precedes Lemma 3.1
\cite{Iv},\cite{Iv-Olsh}.

\begin{lem} \label{lem:betacont}
Let $B \in \mathcal B$ and let $\Gamma$ be a contiguity subdiagram
of a cell $\Pi$ to a $B$-periodic section $p$ of the contour of a
reduced diagram $\Delta$ of rank $i$, such that $lab(\Gamma \wedge
p)$ is a subword of $B^{\pm 3}$. Then $r(\Gamma)=0$ and the
contiguity degree of $\Pi$ to $p$ via $\Gamma$ is less than
$\beta$.
\end{lem}

\begin{proof}
We prove this lemma by contradiction. Assume that the triple
$(p,~\Gamma,~\Pi)$ is a counterexample, where the contiguity
subdiagram $\Gamma$ has a minimal type. Let the standard contour
of $\Gamma$ be $\partial \Gamma =d_1p_1d_2q$, where $p_1=\Gamma
\wedge p$, $q=\Gamma \wedge \Pi$. The bonds defining $\Gamma$ are
$0$-bonds (otherwise $(p,~\Gamma,~\Pi)$ is not minimal). This
means that $~\max\{\|d_1\|,~\|d_2\|\} <\Lambda$. Assuming that
$r(\Gamma)>0$, by Lemma 5.7 \cite{Iv},\cite{Iv-Olsh}
($\tau(\Gamma)<\tau(\Delta)$), there is a $\theta$-cell. However,
by minimality of $(p,~\Gamma,~\Pi)$, the degree of contiguity of
any cell from $\Gamma$ to $p_1$ is less than $\beta$ and the
contiguity degree of any cell from $\Gamma$ to $q$ is less than
$\alpha$ by Lemma 3.4 \cite{Iv},\cite{Iv-Olsh} (again,
$\tau(\Gamma)<\tau(\Delta)$). Contiguity degrees of cells of
positive ranks to sections $d_1$ and $d_2$ of $\partial \Gamma$
are bounded from above by $\beta$ in view of inequality
$~\max\{\|d_1\|,~\|d_2\|\}<\Lambda$ and Lemma 6.1
\cite{Iv},\cite{Iv-Olsh}. The fact that $\theta
> \alpha + 3 \beta$ implies that $\Gamma$ does not have cells of positive ranks, i.e $\Gamma$ is a diagram over the presentation of
$G$. Note that the endpoints of the path $p_1$ are within $\nu_0$
from the geodesic (in the Cayley graph of $G$) segment that
connects the endpoints of $p$. Assuming that the contiguity degree
of $\Pi$ to $p$ via $\Gamma$ is greater than or equal to $\beta$,
we arrive at a contradiction with $(\Lambda+\nu_0, \beta
n)$-aperiodicity of $B^3$ in $G$. Lemma is proved.
%\footnote{tut nado byt' ostorozhnee: dlya polucheniya
%protivorechiya s aperiodichnost'yu nuzhno, chtoby primykanie bylo
%k geodezicheskomu puti. A nash put' $p$, navernoe, mozhet
%geodezicheskim i ne yavlyat'sya. Hotya, po Lemme
%\ref{lem:cyc_min-qd}, on budet
%$(\lambda_0,c_0)$-kvazigeodezichekim (\club v grafe cayley kakoi
%gruppy? $G(i)$?), a znachit, po Lemme \ref{lem:close}, on budet
%$\nu_0$-blizok k geodezicheskoj. T.e. v samom nachale sekcii
%dostatochno potrebovat' chtoby slova $B\in \mathcal B$ byli
%$(\Lambda+\nu_0,\beta n)$-aperiodicheskimi.}
\end{proof}

The formulation of the analogue of Lemma 15.10
\cite{Iv},\cite{Iv-Olsh} is changed as follows. For any finite
subgroup $\mathcal G$ of $G_{\omega}(i)$ equations from
\cite{Iv},\cite{Iv-Olsh} hold with the exponent $n$ replaced by
$2n$. If, in addition, $\mathcal G$ is a finite subgroup of
$G_{\omega}(i)$, conjugate to a subgroup of $\mathcal K(A_j)$, $j
\le i$, where the period $A_j$ is special,
%such that a non-trivial power of
%some word $B_k \in \mathcal B$ is conjugate in $G_{\omega}(j-1)$
%to an element $A_j^{\ell}F$, for some $\ell \in \mathbb Z, \ell
%\ne 0, F \in \mathcal F(A_j) \subseteq G_{\omega}(j-1)$ and
%$\omega_k=0$ (in other words, the rank $j$ relation is $A_j^n=1$),
then the equations from part (a) of Lemma 15.10
\cite{Iv},\cite{Iv-Olsh} hold with exponent $n$. The proof is
retained. In the case the relation of rank $j$ is $A_j^n=1$, the
claim follows from Lemmas \ref{lem:trivsbgp} and \ref{lem:noinvol}
(below) applied in smaller rank: the subgroup $\mathcal F(A_j)$ is
trivial, there are no $\mathcal F(A_j)$-involutions, and, hence,
the equations are trivially satisfied.

%%%%%%%%%%%%%%%%%%%%%%%%%%%%%%%%%%%%%%%%%%%%%%%%%%%%%%%%%%%%%%%%%%%%%%%%%%%%%%%%%%%%%%

The inductive step from rank $i$ to $i+1$ starts with the
following lemmas.

\begin{lem} \label{lem:inforder_simple} If a word $B$ from $\mathcal B$ is of infinite order in rank $i$, then $B$
is simple in rank $i$.
\end{lem}

\begin{proof}
By the choice of the set $\mathcal B$ and the definition of simple
words, $B$ is simple in rank $0$ and therefore it is simple in all
ranks $i \le i_0$ (as in \cite{Iv-Olsh}, $i_0$ stands for the
maximal rank for which $\|A_{i_0}\|<C$). Let now $i
>i_0$. Assume that $B$ is not simple in rank $i$. By definition
of a simple in rank $i$ word, and because $B$ has infinite order
in rank $i$, this means that $B$ is not cyclically reduced in rank
$i$, that is, there exists a word $X$ conjugated to $B$ in rank
$i$ such that $\|X\|<\|B\|$. Let $\Delta$ be an annular reduced
diagram representing conjugacy of $B$ and $X$ in rank $i$. Denote
by $p, ~q$ contours of $\partial \Delta$ labelled by $B,~X$
respectively. Without loss of generality we may assume that the
path $q$ is cyclically geodesic in $\Delta$, i.e. geodesic in its
homotopy class in $\Delta$. If $r(\Delta)>0$ then, by Lemma 5.7
\cite{Iv},\cite{Iv-Olsh}, there is a cell $\Pi$ of positive rank
in $\Delta$ and contiguity subdiagrams $\Gamma_p$, $\Gamma_q$ of
$\Pi$ to $p$, $q$ respectively such that the sum of contiguity
degrees of $\Pi$ to $p, ~q$ via $\Gamma_p,~\Gamma_q$ is greater
than $\theta$. However, by Lemma \ref{lem:betacont}, $\|\Gamma_p
\wedge \Pi\|<\beta \|\partial \Pi\|$ and, by Lemma 3.3
\cite{Iv},\cite{Iv-Olsh}, $\|\Gamma_q \wedge \Pi\|<\alpha
\|\partial \Pi\|$. Thus, the inequality $\theta>\alpha+\beta$
implies that $\Delta$ does not have cells of positive rank. But
this is a contradiction since the word $B$ was chosen to be
cyclically reduced in $G=G(0)$.
\end{proof}

\begin{lem} \label{lem:inforder_noncom}
Let $l \ne k$ and assume that some words $B_l,~B_k \in \mathcal B$
are of infinite order in the group $G_{\omega}(i)$. Then $B_l$ and
$B_k$ are not commensurable in $G_{\omega}(i)$.
\end{lem}

\begin{proof}
Let $l>k$. By the choice of the set $\mathcal B$, the words $B_l$
and $B_k$ are not commensurable in $G=G_{\omega}(0)$. Let now
$i\ge 1$, and assume, on the contrary, that the words $B_l, B_k
\in \mathcal B$ have infinite order and are commensurable in
$G_{\omega}(i)$. Let $\Delta$ be an annular reduced diagram of
rank $i$ for the conjugacy of some non-trivial powers $B_l^s$ and
$B_k^t$. Denote by $p$ and $q$ contours of $\Delta$ so that
$lab(p)\equiv B_l^s$, $lab(q) \equiv B_k^t$. The words $B_l,~B_k$
are not commensurable in $G_{\omega}(0)$, hence $r(\Delta)>0$. By
Lemma~5.7 \cite{Iv},\cite{Iv-Olsh}, there is a cell $\Pi$ of
positive rank $j \le i$ and contiguity subdiagrams $\Gamma_p$,
$\Gamma_q$ of $\pi$ to $p$, $q$ respectively, such that
$\|\Gamma_p \wedge \Pi\|+\|\Gamma_q \wedge \Pi\|>\theta\|\partial
\Pi\|$.

Note that $\|A_j\| \le \|B_k\|<n^{-2}\|B_l\|$. Therefore
$\|\partial \Pi\|<\|B_l\|$ and $lab(\Gamma_p \wedge p)$ is a
subword of $B_l^{\pm 3}$. It follows from Lemma \ref{lem:betacont}
that $\|\Gamma_p \wedge \Pi\|< \beta \|\partial \Pi\|$. By Lemmas
\ref{lem:inforder_simple} and 3.4
\cite{Iv},\cite{Iv-Olsh}
%\footnote{a v usloviyah Lemmy 3.4 ne
%trebuetsya, chtoby put' $q$ gladkim? Esli da, to pochemu jeto
%vyploneno v nashem sluchae? \club potomu chto B - simple, i po
%opredeleniyu, put', metka kotorogo - stepen' prostogo slova,
%yavlyaetsya gladkim}
, $\|\Gamma_q \wedge \Pi\|< \alpha \|\partial
\Pi\|$. This is a contradiction since $\theta > \alpha + \beta$.
\end{proof}

%%%%%%%%%%%%%%%%%%%%%%%%%%%%%%%%

The following lemmas describe the structure of finite subgroups
associated to the periods distinguished by the elements from the
set $\mathcal B$. They precede Lemma 18.5
\cite{Iv},\cite{Iv-Olsh}.

\begin{lem} \label{lem:trivsbgp}
Let a non-trivial power $B^r$ of some word $B \in \mathcal B$ be
conjugate in $G_{\omega}(i)$ to an element $A_{i+1}^{\ell}F$, for
some $\ell \in \mathbb Z, \ell \ne 0, F \in \mathcal F(A_{i+1})$.
Then $\mathcal F(A_{i+1})=\{1\}$.
\end{lem}

\begin{proof}
The word $B$ is of infinite order in rank $i$ and, by Lemma
\ref{lem:inforder_simple} it is simple in rank $i$. Let $K$ be a
finite subgroup of $G_{\omega}(i)$ normalized by $B^r$. Our aim is
to show that $K$ is trivial. Assuming the contrary, choose a word
$F_0$ representing a non-trivial element from $K$. Since $K$ is
finite and is normalized by $B^r$, there exists an integer $s>0$
such that
$$
B^{-s}F_0B^s\eqi F_0.
$$
Choose $s$ such that $|s| > n^2\|F_0\|$ and let $\Delta$ be a disk
reduced diagram of rank $i$ representing this equation with the
standard contour $\partial \Delta=bpcq$, where $lab(b)\equiv
F_0,lab(p) \equiv B^s, lab(c) \equiv F_0^{-1}, lab(q) \equiv
B^{-s}$. Assuming $r(\Delta)>0$, by Lemma
5.7\cite{Iv},\cite{Iv-Olsh}, one will find a $\theta$-cell $\Pi$
in $\Delta$. Denote $j=r(\Pi)$ (remark that $\|B\|\ge \|A_j\|$ and
let $\Gamma_b, \Gamma_p, \Gamma_c$ and $\Gamma_q$ be contiguity
subdiagrams (some of them may be absent) of $\Pi$ to $b,p,c$ and
$q$ respectively. Sections $b$ and $c$ may be regarded as geodesic
sections of $\partial \Delta$. Therefore, by Lemma 3.3
\cite{Iv},\cite{Iv-Olsh},
$$\max\{\|\Gamma_b \wedge \Pi\|,~\|\Gamma_c \wedge \Pi\|\}<\alpha
\|\partial \Pi\|.$$

Suppose that $\|\Gamma_b \wedge \Pi\|+\|\Gamma_c \wedge \Pi\| >
(\alpha+\beta)\|\partial \Pi\|$. In particular, both $\Gamma_b$
and $\Gamma_c$ are present, and the contiguity degrees of $\Pi$ to
$b$ and to $c$ via $\Gamma_b$ and $\Gamma_c$ are greater than
$\beta$. Denote the standard contours $\partial
\Gamma_b=u_1b_1u_2u_b$ and $\partial \Gamma_c=v_1c_1v_2u_c$, where
$b_1 = \Gamma_b \wedge b$, $u_b = \Gamma_b \wedge \Pi$, $c_1 =
\Gamma_c \wedge c$, $u_c = \Gamma_c \wedge \Pi$. By Lemma 3.1
\cite{Iv},\cite{Iv-Olsh},
$$\max\{\|u_1\|, \|u_2\|, \|v_1\|, \|v_2\|\}<\gamma \|\partial \Pi\|.$$ It
follows that initial and terminal vertices $p_-$ and $p_+$ of the
path $p$ can be joined in $\Delta$ by a path of length less than
$$
\|F_0\|+\gamma \|\partial \Pi\|+1/2\|\partial \Pi\|+\gamma
\|\partial \Pi\|+\|F_0\|=2\|F_0\|+(2\gamma+1/2)n\|A_j\|<\rho
\|p\|.
$$
This contradicts Lemma 6.1 \cite{Iv},\cite{Iv-Olsh} since $p$ is a
smooth section of $ \partial \Delta$ and of the boundary of any
subdiagram of $\Delta$ that contains $p$ in its contour.
Consequently, $\|\Gamma_b \wedge \Pi\|+\|\Gamma_c \wedge \Pi\| \le
(\alpha+\beta)\|\partial \Pi\|$.

Let us assume now that $\Gamma_p$ is present. Denote the standard
contour $\partial \Gamma_p=d_1p_1d_2p_2$, where $p_1= \Gamma_p
\wedge p$ and $p_2=\Gamma_p \wedge \Pi$. If $\|p_2\|\ge \beta
\|\partial \Pi\|$, then by Lemma~\ref{lem:inforder_simple} and
Lemma 20.2 \cite{Iv},\cite{Iv-Olsh} applied to $\Gamma_p$
($\tau(\Gamma_p)<\tau(\Delta)$), one gets
$\|p_1\|<(1+\varepsilon)\|B\|$. Hence, $lab(p_1)$ is a subword of
$B^{\pm 3}$, and Lemma~\ref{lem:betacont} yields a contradiction.
Therefore, $\|\Gamma_p \wedge \Pi\|< \beta \|\partial \Pi\|$.
Arguing in the same way, we obtain that $\|\Gamma_q \wedge \Pi\|<
\beta \|\partial \Pi\|$ in the case $\Gamma_q$ is present. It
follows that $$\|\Gamma_b \wedge \Pi\|+\|\Gamma_p \wedge
\Pi\|+\|\Gamma_c \wedge \Pi\|+\|\Gamma_q \wedge \Pi\|< (\alpha +
3\beta) \|\partial \Pi\| <\theta \|\partial \Pi\|,$$ contradicting
with the definition of a $\theta$-cell.

Therefore, $r(\Delta)=0$ and equality $B^{-s}F_0B^s=F_0$ holds in
the group $G(0)$. This means that $F_0$ belongs to the elementary
subgroup of $B$ in $G(0)$. Since the elementary subgroup of $B$ is
cyclic, $F_0$ is trivial in $G_{\omega}(i)$. Consequently,
$\mathcal F(A_{i+1})=\{1\}$.
\end{proof}

\begin{lem} \label{lem:noinvol}
Let a non-trivial power $B^r$ of some word $B \in \mathcal B$ be
conjugate in $G_{\omega}(i)$ to an element $A_{i+1}^{\ell}F$, for
some $\ell \in \mathbb Z, \ell \ne 0, F \in \mathcal F(A_{i+1})$.
Then there are no $\mathcal F(A_{i+1})$-involutions.
\end{lem}

\begin{proof}
By Lemma \ref{lem:trivsbgp}, we may assume that $B^r \eqi
X^{-1}A_{i+1}^{\ell}X$ for some word $X$. Let $J$ be a $\mathcal
F(A_{i+1})$-involution. Then $J_0^{-1}B^rJ_0 \eqi B^{-r},$ where
$J_0 \eqi X^{-1}JX$. It follows from Lemma \ref{lem:trivsbgp} and
the definition of a $\mathcal F(A_{i+1})$-involution that $J$ (and
therefore $J_0$) is of order two in rank $i$.

Choose $s$ such that $|s|>n^2\|J_0\|$ and consider a reduced
diagram $\Delta$ of rank $i$ representing the equation
$$J_0^{-1}B^sJ_0 \eqi B^{-s}.$$
Arguing as in the proof of Lemma \ref{lem:trivsbgp}, we obtain
that $r(\Delta)=0$. Consequently, $J_0^{-1}B^sJ_0 \eqin B^{-s}$,
and $J_0$ belongs to the elementary subgroup of $B$ in $G(0)$.
Since the order of $J_0$ in $G_{\omega}(i)$ is finite while the
order of $B$ in $G_{\omega}(i)$ is infinite, it follows that
$J_0\eqi 1$ (and therefore $J \eqi 1$).
\end{proof}

The inductive proof is completed as in \cite{Iv}, \cite{Iv-Olsh}.
As a result, we obtain a presentation
$$
G_{\omega}(\infty)=\langle \A~\|~R \in \mathcal R_0 \cup
\{A_1^{n(1)},A_2^{n(2)}, \dots~\}\rangle.
$$
of an infinite group of exponent $2n$. The proof that
$G_{\omega}(\infty)$ has trivial center repeats the argument from
\S~21 \cite{Iv}.

Now we explain that there is a continuum of pairwise
non-isomorphic groups among $\{G_{\omega}(\infty)\}_{\omega}$. We
notice first that, by Lemmas 10.1 and 10.2 \cite{Iv-Olsh},
\cite{Iv}, and Lemma \ref{lem:inforder_noncom}, the set of special
periods is infinite and there is a special period corresponding to
each element $B \in \mathcal B$. Presentations of groups
$G_{\omega}(\infty)$ are now distinguished by the orders of
special periods: if $\omega$ and $\omega'$ are two different
sequences of $0$'s and $1$'s, find the first place (say, $k$)
where they differ. By construction and Lemma 10.4 \cite{Iv-Olsh},
\cite{Iv}, the special period that corresponds to $B_k$ has
different orders in $G_{\omega}(\infty)$ and
$G_{\omega'}(\infty)$, meaning that $G_{\omega}(\infty)$ and
$G_{\omega'}(\infty)$ are quotients of $G$ over different normal
subgroups. By Remark \ref{rem:cont_quot_non-isom}, we see that
there is a continuous family of pairwise non-isomorphic infinite
centerless quotients of $G$ of exponent $2n$. Theorem
\ref{thm:contin_3} is proved.

%%%%%%%%%%%%%%%%%%%%%%%%%%%%%%%%%%%%%%%%%%%%%%%%%%%%%%%%%%%%%%%%%%%%%%%%%%%%%%%%
\section{Periodic quotients of large groups}\label{sec:large}
This section is devoted to proving Theorem \ref{thm:quot-large}
and Corollaries
\ref{cor:large_quot_explicit},\ref{cor:more_gens_than_rels}.

Let $\mathcal C$ be a collection of groups. A group $G$ is said to
be {\it residually in} $\mathcal C$ if for every $g \in
G\setminus\{1\}$ there is a homomorphism $\varphi:G \to H$ for
some group $H \in \mathcal{C}$ such that $\varphi(g) \neq 1$. For
a prime number $p$, denote by $\mathcal{C}_p$ the collection
consisting of all finite $p$-groups. Since every finite $p$-group
is nilpotent we can make

\begin{rem}\label{rem:delta-series-rf} 1) If $H$ is a finite $p$-group, then there is an integer $n \in \N$ such that
$\delta^p_n(H)=\{1\}$. 2) For a finitely generated group $G$,
$\bigcap_{j=0}^\infty \delta^p_j (G)=\{1\}$ if and only if $G$ is
residually in
$\mathcal{C}_p$. %(finite $p$-groups).
\end{rem}

Since free groups are residually in $\mathcal{C}_p$ for any prime
$p$ (\cite[14.2.2]{K-M}), we have
\begin{lem} \label{lem:F-p-approx} If $F$ is a free group
and $p \in \N$ is a prime number then ~~$\displaystyle
\bigcap_{j=0}^\infty \delta^p_j (F)=\{1\}$.
\end{lem}

Throughout this section $F_2$ will denote the free group of rank
$2$. For a group $G$, we will write $G \twoheadrightarrow F_2$ if
there exists an epimorphism from $G$ to $F_2$. The lemma below was
proved by Olshanskii and Osin in \cite{Olsh-Osin}. (We apply it in
the special case when $\Pi=(p,p,\dots)$ is a constant sequence.)

\begin{lem} {\rm (\cite[Lemma 3.2]{Olsh-Osin})}\label{lem:large_delta}
Let G be a finitely generated group, $N \lhd G$ be a normal
subgroup of finite index and let $p \in \N$ be a prime number.
Suppose that $N \twoheadrightarrow F_2$. Then for any element $g
\in N$ there is $m >0 $ such that if $g^n \in \delta^p_m(N)$, then
$\delta^p_m(N)/\langle\langle g^n \rangle\rangle^G
\twoheadrightarrow F_2$.
\end{lem}

(Here we modified the original formulation of \cite[Lemma
3.2]{Olsh-Osin} by replacing $\delta^p_r(G)$ with an arbitrary
finite index normal subgroup $N$, observing that the proof
continues to be valid in this, more general, situation).

\begin{lem} \label{lem:quot-large} Suppose that $G$ is a finitely generated group, $N \lhd G$ is a finite index normal subgroup,
such that $N \twoheadrightarrow F_2$, and $p$ is an arbitrary
prime number. Then for each sequence $\omega \in \Omega = \{0,1\}^\N$, there is
a group $H_\omega$ containing a normal subgroup $Q_\omega \lhd
H_\omega$ such that
\begin{itemize}
\item[(i)] $H_\omega$ is a quotient of $G$; \item[(ii)]
$H_\omega/Q_\omega \cong G/N$; \item[(iii)] $Q_\omega$ is a
periodic $p$-group; \item[(iv)] $\displaystyle
\bigcap_{i=0}^\infty \delta^p_i (Q_\omega)=\{1\}$; \item[(v)] if
$\omega\neq \omega' \in \Omega$ then there is $ v \in \N$ such
that $|Q_\omega/\delta^p_v(Q_\omega)|\neq
|Q_{\omega'}/\delta^p_v(Q_{\omega'})|$; consequently $Q_\omega
\not\cong Q_{\omega'}$.
\end{itemize}
\end{lem}

\begin{proof} Our argument will be similar to the one used in the proof of \cite[Thm. 1.2]{Olsh-Osin}.

First, enumerate all the elements of $N$: $N=\{f_1,f_2,\dots \}$.
Let $\Omega_i=\{0,1\}^i$ be the set of sequences of $0$'s and
$1$'s of length $i$, $i \in \N$, and let $\Omega_0$ consist of the
empty sequence $\emptyset$. Denote $G_\emptyset=G$, $r_\emptyset
=q_\emptyset=0$, and assume that for some $i \ge 0$ and for every
$\iota \in \Omega_i$ we have already constructed the
quotient-group $G_\iota$ of $G$ and integers $r_\iota, q_\iota \ge
0$ such that the images of $f_1,\dots,f_i$ have finite orders in
$G_\iota$ and $\delta^p_{r_\iota}(N_\iota) \twoheadrightarrow
F_2$, where $N_\iota$ is the image of $N$ in $G_\iota$ under the
natural homomorphism (the numbers $q_\iota$ are auxiliary, and
will be used during the inductive argument).

Choose any $\zeta=(\zeta_1,\dots,\zeta_{i+1}) \in \Omega_{i+1}$,
and take $\iota=(\zeta_1,\dots,\zeta_i) \in \Omega_i$.

Observe that the subgroup $E_\iota=\bigcap_{j=0}^\infty
\delta^p_j(N_\iota)$ is characteristic in $N_\iota$, hence it is
normal in $G_\iota$. Denote $K_\iota=G_\iota/E_\iota$,
$L_\iota=N_\iota/E_\iota \lhd K_\iota$; then
\begin{equation} \label{eq:L-iota} \bigcap_{j=0}^\infty \delta^p_j(L_\iota)=\{1\} .\end{equation}
Since $\displaystyle E_\iota=\bigcap_{j=0}^\infty
\delta^p_j\left(\delta^p_{r_\iota}(N_\iota)\right)$, Lemma
\ref{lem:F-p-approx} implies that $\delta^p_{r_\iota}(L_\iota)
\twoheadrightarrow F_2$. Let $l \ge i+1$ be the smallest index
such that the image of $f_l$ has infinite order in $K_\iota$.
Denote by $g$ the image of $(f_l)^{p^{r_\iota}}$ in $K_\iota$;
thus $g \in \delta^p_{r_\iota}(L_\iota)$ is an element of infinite
order.

Now we apply Lemma \ref{lem:large_delta} to find $m \in \N$ such
that
\begin{equation} \label{eq:large-quot}\delta^p_{r_\iota+m}(L_\iota)/\langle \langle g^n \rangle \rangle^{K_\iota}=
\delta^p_{m}\left( \delta^p_{r_\iota}(L_\iota) \right)/\langle
\langle g^n \rangle \rangle^{K_\iota} \twoheadrightarrow
F_2,\end{equation} for any integer $n$ such that $g^n \in
\delta^p_{r_\iota+m}(L_\iota)$.

Set $s=\max\{r_\iota+m,q_\iota\}$ so that $g^{p^s} \in
\delta^p_{s}(L_\iota) \le \delta^p_{r_\iota+m}(L_\iota) \cap
\delta^p_{q_\iota}(L_\iota)$. By \eqref{eq:L-iota} one can find
the smallest integer $v_\iota$ such that $g^{p^s} \notin
\delta^p_{v_\iota}(L_\iota)$. Set $r_{\zeta}=r_\iota+m$ and
$q_\zeta=v_\iota$. Note that $q_\zeta>s\ge q_\iota$, and the
integers $r_\zeta$ and $q_\zeta$ only depend on the prefix $\iota$
of the sequence $\zeta$, and are independent of the value of
$\zeta_{i+1}$.

If $\zeta_{i+1}=0$, define $n_\zeta=p^s$ (thus, $g^{n_\zeta} \in
L_\iota \setminus \delta^p_{v_\iota}(L_\iota)$), otherwise, if
$\zeta_{i+1}=1$, define ${n_\zeta}=p^{q_\zeta}$ (thus,
$g^{n_\zeta} \in \delta^p_{v_\iota}(L_\iota)$). Denote
$G_\zeta=K_\iota/\langle \langle g^{n_\zeta} \rangle
\rangle^{K_\iota}$, and let $N_\zeta$ be the image of $L_\iota$ in
$G_\zeta$ (see the diagram below).

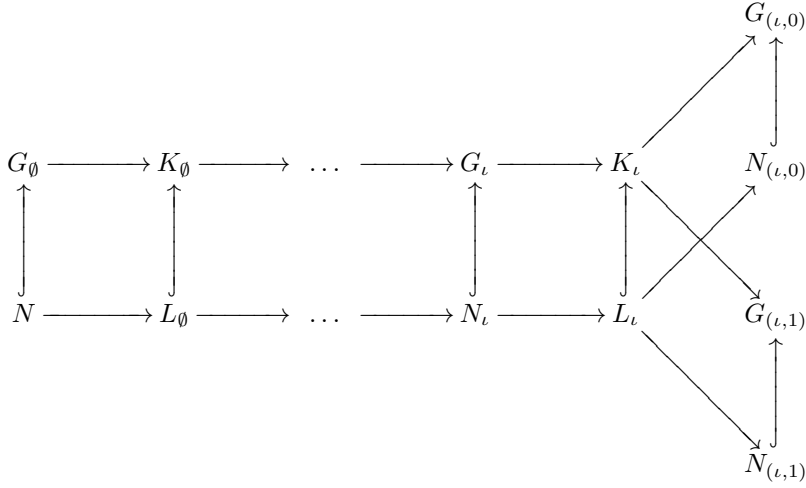
\Diagram
              & & & & & & & & & &G_{(\iota,0)}\\
              & & & & & & & & &\ruTo &\uInto \\
G_\emptyset &\rTo &K_\emptyset &\rTo &\mbox{ }\dots\mbox{ } &\rTo &G_\iota & \rTo & K_\iota & &N_{(\iota,0)}\\
\uInto & & \uInto & & & &\uInto & &\uInto &\rdTo\ruTo & \\
N &\rTo &L_\emptyset &\rTo &\mbox{ }\dots\mbox{ } &\rTo &N_\iota & \rTo & L_\iota & &G_{(\iota,1)}\\
              & & & & & & & & &\rdTo &\uInto \\
              & & & & & & & & & &N_{(\iota,1)}\\
\endDiagram

Observe, that
$$\delta^p_{r_\zeta}(N_\zeta)\twoheadrightarrow F_2 ~\mbox{ by \eqref{eq:large-quot}},$$
\begin{equation}\label{eq:same_quot} N_\zeta/\delta^p_{u}(N_\zeta) \cong L_\iota/\delta^p_{u}(L_\iota)
\cong N_\iota/\delta^p_{u}(N_\iota) ~\mbox{ for every } u \le
q_\iota ,~\mbox{ and} \end{equation}
% $$ |N_\zeta/\delta^p_{v_\iota}(N_\zeta)| = |N_\iota/\delta^p_{v_\iota}(N_\iota)|
%~\mbox{ if and only if } \zeta_{i+1}=1, ~\mbox{ because} $$
\begin{equation}\label{eq:dist_quot} |N_{(\iota,0)}/\delta^p_{v_\iota}(N_{(\iota,0)})|<|N_{(\iota,1)}/\delta^p_{v_\iota}(N_{(\iota,1)})|=
|L_{\iota}/\delta^p_{v_\iota}(L_{\iota})|=|N_{\iota}/\delta^p_{v_\iota}(N_{\iota})|<\infty.\end{equation}

Thus, for each $i \ge 0$ and each $\iota \in \Omega_i$, we have
constructed a quotient $G_\iota$ of $G$. Additionally, we have an
epimorphism $G_\iota \twoheadrightarrow G_\zeta$ whenever $\iota$
is a prefix of $\zeta$. Therefore, for every $\omega \in \Omega$
we can define the group $G_\omega$ as a direct limit of the groups
$G_\iota$, where $\iota$ runs over all prefixes of $\omega$. Let
$N_\omega \lhd G_\omega$ be the image of $N$ in $G_\omega$, and
let $H_\omega$ be the quotient of $G_\omega$ by
$E_\omega=\bigcap_{j=0}^\infty \delta^p_j(N_\omega)$. Finally, let
$Q_\omega$ be the image of $N_\omega$ in $H_\omega$.

By construction, $Q_\omega$ is periodic (as a quotient of $N$),
and $\bigcap_{j=0}^\infty \delta^p_j(Q_\omega)=\{1\}$. Since $N$
is finitely generated, then so is $Q_\omega$. Therefore, for every
$j$, $Q_\omega/\delta^p_j(Q_\omega)$ is a finite group whose order
divides $p^j$, hence $Q_\omega$ is a $p$-group.

The second claim of the lemma also follows from the construction.
It remains for us to prove the fifth claim. Suppose $\omega \neq
\omega' \in \Omega$, and let $\iota \in \Omega_i$, $i \ge 0$, be
the longest common prefix of $\omega$ and $\omega'$.
%Arguing by contradiction, assume that $Q_\omega \cong Q_{\omega'}$. Consequently
%\begin{equation}\label{eq:for_contrad}|Q_\omega/\delta^p_{v_\iota}(Q_\omega)| = |Q_{\omega'}/ \delta^p_{v_\iota}(Q_{\omega'})|.
%\end{equation}

By \eqref{eq:same_quot} we have
$$Q_\omega/\delta^p_{v_\iota}(Q_\omega) \cong N_\omega/\delta^p_{v_\iota}(N_\omega) \cong N_\zeta/\delta^p_{v_\iota}(N_\zeta)~\mbox{ and }$$
$$Q_{\omega'}/\delta^p_{v_\iota}(Q_{\omega'}) \cong N_{\omega'}/\delta^p_{v_\iota}(N_{\omega'})
\cong N_{\zeta'}/\delta^p_{v_\iota}(N_{\zeta'}), $$ where
$\zeta,\zeta' \in \Omega_{i+1}$ are prefixes of length $i+1$ of
$\omega$ and $\omega'$ respectively. Using this together with
\eqref{eq:dist_quot} we get
$$|Q_\omega/\delta^p_{v_\iota}(Q_\omega)|=|N_\zeta/\delta^p_{v_\iota}(N_\zeta)| \neq |N_{\zeta'}/\delta^p_{v_\iota}(N_{\zeta'})|=
|Q_{\omega'}/ \delta^p_{v_\iota}(Q_{\omega'})|,$$ yielding (v).
Thus the lemma is proved.
\end{proof}

\begin{proof}[Proof of Theorem \ref{thm:quot-large}] %First, choose an arbitrary $p \in P$.
Applying Lemma \ref{lem:quot-large} to $G$, $N$ and $p$, for each
$\omega \in \Omega$ we obtain a group $H_\omega$ together with a
normal subgroup $Q_\omega$ of index $k=|G:N|$ enjoying the
properties (i)-(v) from the claim of Lemma~\ref{lem:quot-large}.

Since $Q_\omega$ is finitely generated, property (iv) implies that
it is residually finite. And since
$|H_\omega:Q_\omega|=|G:N|<\infty$ by (ii), $H_\omega$ is
residually finite for every $\omega \in \Omega$. Property (iii)
ensures that $Q_\omega$ is a $p$-group, and, hence, the group
$H_\omega$ is periodic.

It remains to observe that for every $\omega \in \Omega$, the
finitely generated group $G_\omega$ can have only finitely many
distinct normal subgroups of index $k$, hence among
$\{G_\omega~|~\omega \in \Omega\}$ there must be
$|\Omega|=2^{\aleph_0}$ of pairwise non-isomorphic groups.
\end{proof}

\begin{proof}[Proof of Corollary \ref{cor:large_quot_explicit}]
First, we observe that since $G/N$ is a $p$-group, there is $n\in
\N$ such that $\delta_n^p(G) \le N$ (by the first part of Remark
\ref{rem:delta-series-rf}). Since $|G:\delta_n^p(G)|<\infty$ and
there is an epimorphism $\varphi:N \twoheadrightarrow F_2$, the
subgroup $\varphi(\delta_n^p(G))$ has finite index in
$\varphi(N)=F_2$. Hence $\delta_n^p(G)\twoheadrightarrow F_2$.
Thus, without loss of generality, we can suppose that
$N=\delta_n^p(G)$ for some $n\in \N$.

Now, we apply Lemma \ref{lem:quot-large} to $G$, $N=\delta_n^p(G)$
and $p$, to obtain the groups $H_\omega$ as in its claim. Since
$H_\omega/Q_\omega\cong G/\delta_n^p(G)$ is a finite $p$-group,
the property (iii) of $Q_\omega$ implies that $H_\omega$ is a
periodic $p$-group. From the proof of Lemma \ref{lem:quot-large}
we see that $Q_\omega$ is the image of $N$ under the natural
homomorphism $G \to H_\omega$, therefore
$Q_\omega=\delta^p_n(H_\omega)$, and for every $j\in \N$,
$\delta^p_j(Q_\omega)=\delta^p_{j+n}(H_\omega)$. Hence, property
(iv) of $Q_\omega$ yields property (b) of $H_\omega$, and
properties (ii) and (v) together yield (c) (for $t=v+n$).
\end{proof}

\begin{proof}[Proof of Corollary \ref{cor:more_gens_than_rels}] By a theorem of Baumslag and Pride~\cite{Baumslag-Pride}, for any prime $p$
there are a normal subgroup $N \lhd G$ and $m \in \N$ such that
$|G:N|=p^m$ and $N \twoheadrightarrow F_2$. To conclude, it
remains to apply Corollary~\ref{cor:large_quot_explicit} to the
pair $(G,N)$.
\end{proof}

%%%%%%%%%%%%%%%%%%%%%%%%%%%%%%%%%%%%%%%%%%%%%%%%%%%%%%%%%%%%%%%%%%%%%%%%%%%%%%%%%%%

{\bf Acknowledgements.} The first and the third authors are grateful to
Institut des Hautes \'Etudes Scientifiques (IH\'ES) for
the hospitality during their work on this paper. The authors
would also like to thank L. Bartholdi and R.I. Grigorchuk for interesting discussions, and the referee for his remarks.

%%%%%%%%%%%%%%%%%%%%%%%%%%%%%%%%%%%%%%%%%%%%%%%%%%%%%%%%%%%%%%%%%%%%%%%%%%%%%%%%%%%

\end{document}